\documentclass[a4paper,12pt,thmsa]{article}

\usepackage{amsfonts}
\usepackage{amssymb}
\usepackage{graphicx}
\usepackage{amsmath}
\usepackage{xr}
\usepackage{bezier}
\usepackage{tikz}
\usepackage[normalem]{ulem}

\makeatletter

\renewcommand{\paragraph}{\@startsection{paragraph}{4}{0mm}{-3mm}{-3mm} {\noindent \bf}}

\renewcommand{\subparagraph}{\@startsection{paragraph}{4}{3mm}{-3mm}{-3mm} {\noindent \bf}}

\newcommand{\RR}{{\mathbb{R}}}

\newcommand{\CF}{{\mathcal{F}}}
\newcommand{\CG}{{\mathcal{G}}}

\newcommand{\CL}{{\mathcal{L}}}
\newcommand{\CT}{{\mathcal{T}}}
\newcommand{\CB}{{\mathcal{B}}}

\newtheorem{thm}{Theorem}
\newtheorem{cor}{Corollary}
\newtheorem{lem}{Lemma}
\newtheorem{pro}{Proposition}

\let\inf\relax \DeclareMathOperator*\inf{\vphantom{p}inf}

\begin{document}

\hyphenchar\font=-1

\baselineskip=7mm

\vspace{3cm}

\title{\textbf{Investment Timing and Technological Breakthroughs}\thanks{We thank St\'ephane Auray, Giorgio Ferrari, Frank Riedel, and Jan-Henrik Steg for extremely valuable feedback. We also thank conference participants at the 2019 Universit\"at Bielefeld ZiF Workshop for useful discussions. This research has benefited from financial support of the ANR (Programme d'Investissement d'Avenir ANR-17-EURE-0010), the research foundation TSE-Partnership, and the AI Interdisciplinary Institute ANITI, which is funded by the French "Investing for the
Future - PIA3" program under the Grant agreement ANR-19-PI3A-0004.}}
\author{Jean-Paul D\'ecamps\thanks{Toulouse School of Economics, University of Toulouse Capitole, Toulouse, France. E-mail: \tt{jean-paul.}  \tt{decamps@tse-fr.eu}.}~~~~~~~Fabien Gensbittel\thanks{Toulouse School of Economics, University of Toulouse Capitole, Toulouse, France. E-mail: \tt{fabien.} \tt{gensbittel@tse-fr.eu}.}~~~~~~~Thomas Mariotti\thanks{Toulouse School of Economics, CNRS, University of Toulouse Capitole, Toulouse, France, CEPR, and CESifo. Email: \tt{thomas.mariotti@tse-fr.eu}.}} \vspace{1cm}

\vspace{1cm}

\maketitle \vspace*{6mm}

\begin{abstract}

We study the optimal investment policy of a firm facing both technological and cash-flow uncertainty. At any point in time, the firm can decide to invest in a stand-alone technology or to wait for a technological breakthrough. Breakthroughs occur when market conditions become favorable enough, exceeding a certain threshold value that is ex-ante unknown to the firm. A microfoundation for this assumption is that a breakthrough occurs when the share of the surplus from the new technology accruing to its developer is high enough to cover her privately observed cost. We show that the relevant Markov state variables for the firm's optimal investment policy are the current market conditions and their current historic maximum, and that the firm optimally invests in the stand-alone technology only when market conditions deteriorate enough after reaching a maximum. Empirically,  investments in new technologies requiring the active cooperation of developers should thus take place in booms, whereas investments in state-of-the-art technologies should take place in busts. Moreover, the required return for investing in the stand-alone technology is always higher than if this were the only available technology and can take arbitrarily large values following certain histories. Finally, a decrease in development costs, or an increase in the value of the new technology, makes the firm more prone to bear downside risk and to delay investment in the stand-alone technology.

\bigskip



\end{abstract}

\thispagestyle{empty}
\newpage \setcounter{page}{1}

\section{Introduction}

It has long been emphasized that firms' incentives to adopt a technology depend on their expectations about future technological improvements (Schumpeter (1950)). Thus Rosenberg (1976, page 525) argues that ``as soon as we accept the perspective of the on-going nature of much technological change, the optimal timing of an innovation becomes heavily influenced by expectations concerning the timing and significance of \textit{future} improvements.'' As a result, a firm expecting imminent technological improvements may be reluctant to adopt the current state-of-the-art technology, for fear of committing itself to practices soon to be antiquated. Expectations of future improvements may thus be a factor accounting for the apparently slow rate at which, historically, some technologies---such as the steam engine in shipping, the oxygen steel-making process, or, more recently, photovoltaic power and hybrid vehicles---have diffused throughout the economy (Rosenberg (1972)).

Since the seminal contribution of Balcer and Lippman (1984), models of technology adoption under technological uncertainty have assumed that the process through which technological innovations are made available to firms is exogenous to the industry under consideration. In particular, it is independent of market conditions, such as the evolution of demand for the final goods manufactured by firms in the client industry. This implies that the interactions between the supply and the demand for innovations are not taken into account, even in a reduced-form way. This assumption may be relatively harmless in the case of general-purpose technologies, but less so for technologies that are targeted at a specific sector. For instance, the incentives of aircraft manufacturers to develop new airplanes can hardly be divorced from the evolution of air traffic, and thus from airlines' willingness to invest in new equipment. Or, to take another example, the incentives of arms manufacturers to develop new weapon systems depend on geopolitical factors that affect the DoD's willingness to acquire such weapons.

To speak to these issues, this paper characterizes the optimal investment policy of a firm who can invest in the current, \textit{stand-alone} technology, or wait until a \textit{breakthrough} occurs and a superior technology becomes available. We build on a standard real-options model in the spirit of McDonald and Siegel (1986) or Dixit and Pindyck (1994), in which market conditions---representing, for instance, the output price in the industry under consideration--- evolve in continuous time according to a diffusion process. The novelty of our approach, compared to real-options models of investment under technological and cash-flow uncertainty, is that breakthroughs do not occur independently of market conditions. Specifically, we assume that a breakthrough occurs only when market conditions become favorable enough, exceeding a certain threshold value. Technological uncertainty is captured by assuming that the firm has incomplete information about the value of this threshold. Thus it does not know when, and under which market conditions, a breakthrough will occur.

Because they are key to our results, we provide a microfoundation for these assumptions in Section \ref{AMotExa}. The idea of this motivating example is to explicitly model the supply of the superior technology by competing \textit{developers}. These developers are assumed to be caught in a race to be the first to supply this technology to the firm; supplying the technology involves a development cost, the same for each developer, which is their private information. Following the usual logic of rent dissipation in winner-take-all preemption games, as in Fudenberg and Tirole (1985), each developer is ready, in equilibrium, to develop the new technology as soon as the payoff from doing so covers the development cost. Assuming that the surplus of investing in the superior technology relative to investing in the stand-alone technology is shared between the firm and the developer who wins the race according to the Nash (1950) bargaining solution, this implies that a breakthrough must occur in a boom---specifically, when market conditions reach a certain threshold value, which triggers the development of the new technology. However, because the firm does not know the cost of developing the new technology, it does not know under which market conditions the developers just break even either. This gives the firm's optimal stopping problem the general structure that is analyzed in this paper.

In this context, we study the optimal investment decision of the firm. At each point in time, the firm can invest in the stand-alone technology, or wait until the superior technology becomes available. A key observation is that the resulting optimal stopping problem is not Markovian with respect to current market conditions, unlike in the real-options literature on technological adoption that assumes a memoryless process for technological breakthroughs. This reflects that, at each point in time, the history of market conditions encodes both the desirability of investing in the stand-alone technology, which typically is high when current market conditions are favorable, and what the firm has learned about the likelihood of benefiting from the superior technology. An implication of this is that the optimal investment policy is necessarily path-dependent.

An intuitive property of the optimal investment policy is that it is not optimal for the firm to immediately invest in the stand-alone technology once market conditions reach a new maximum, because this may be precisely the time at which a breakthrough occurs and the superior technology becomes available. By contrast, once such a maximum has been reached and market conditions start to deteriorate, then the firm for a while no longer learns about the advent of a breakthrough---for instance, in our motivating example, it no longer learns about the cost of developing the new technology. This suggests that the relevant Markov state variables are the current market conditions and their current historic maximum, as established in Proposition \ref{markovian}. This result further shows that the option value of waiting for a breakthrough can be written as the expectation of an integral with respect to future increments of the maximum process, which is a novel feature of the optimal stopping problem studied in this paper. An intuitive complement to this result is that, after market conditions have reached a maximum yet no breakthrough has occurred, the firm revises its beliefs and becomes more pessimistic about benefiting from the superior technology, at least in the near future. In particular, the more market conditions deteriorate, the longer it will take for the market to recover and for the firm to resume learning. Our central result, Theorem \ref{maintheo}, shows that, as a consequence, investment in the stand-alone technology takes place when market conditions cross from above an investment boundary.

The resulting optimal investment policy is characterized by two phases. In a first phase, which takes place before the market conditions have reached a certain level, the firm never invests in the stand-alone technology; it may, however, invest in the superior technology, should it become available. Once this level has been reached, the firm enters into a new phase, and invests in the stand-alone technology when market conditions deteriorate enough after having reached a maximum. The corresponding lower investment threshold is not a constant, unlike in Dixit's (1989) model of entry and exit decisions under uncertainty. Rather, this threshold, which is always greater than the optimal stand-alone investment threshold, is a strictly increasing function of the maximum market conditions achieved so far, which describes the optimal investment boundary. This implies that the required return for investing in the stand-alone technology is always higher than if this were the only available technology, and that it can take arbitrarily large values following certain histories. Theorem \ref{maintheo} gives a complete characterization of the optimal investment boundary and of the firm's optimal value function as the solutions to a variational problem; specifically, the optimal investment boundary is shown to be the unique solution of an ordinary differential equation that satisfies a limit condition, akin to a transversality condition, and the firm's value function is explicit given the investment boundary.

In the context of our motivating example, a testable implication of our model is that investments in new technologies requiring the active cooperation of developers should take place in booms, whereas investments in the current state-of-the-art technology should take place in busts. This is because favorable market conditions give developers an incentive to supply superior technologies by increasing the surplus from investing in such technologies rather than in the state-of-the-art technology, which they share with the firms in the client industry. This, in turn, gives incentive to these firms to invest in the state-of-the-art technology only if market conditions deteriorate enough, while still remaining sufficiently favorable so that investment would take place immediately if firms were not expecting future technological breakthroughs. We show that a decrease in development costs in the hazard-rate order, or an increase in the value of the new technology in both absolute and marginal terms, makes the firm more prone to bear downside risk and to delay investment in the stand-alone technology.

\vskip 3mm

The paper is organized as follows. Section \ref{AMotExa} provides our motivating example. Section \ref{RelLit} discusses the relevant economic and mathematical literature. Section \ref{TheMod} precisely describes the firm's problem and the assumptions under which we solve it. Section \ref{AMaFor} provides the Markovian formulation of the firm's problem. Section \ref{TheMaiThe} heuristically derives and formally states our main theorem. Section \ref{analytical} is devoted to the analysis of the variational system. Section \ref{section_verification} completes the proof of our main theorem by providing the required verification argument. Section \ref{Discussion} discusses the implications of our analysis. Section \ref{Concluding Remarks} concludes. Proofs not given in the main text are collected in Appendices A--C.

\section{A Motivating Example} \label{AMotExa}

Our motivating example is an investment-timing game with technological breakthroughs that can be informally described as follows. The benchmark is that of a decision maker (DM) who decides when to invest in a project whose net value at investment time $\tau$ is $R(X_\tau)$, where $X \equiv (X_t)_{t\geq 0}$ is an observable and continuous Markov process and $R$ is a net payoff function. The interpretation of $X$ depends on the intended application of the model. If the DM is a firm, $X$ may stand for the cash-flow upon investing or for the consumers' willingness-to-pay for the DM's output. If $X$ is a social planner, $X$ stands for the social desirability of the investment project. Whatever the interpretation, the DM, if left to his own devices, would solve the following standard optimal stopping problem:
\begin{eqnarray} \label{bof}
V_R(x) \equiv \sup_{\tau } \, \mathbf E_x \hskip 0.3mm[\mathrm e^{-r\tau } R(X_\tau)],
\end{eqnarray}
where $r>0$ is the DM's discount rate. Now, suppose that the DM can benefit from the help of two competing developers D1 and D2, each of whom can develop, at some endogenously determined random time, a new technology allowing the DM to improve his net payoff function from $R$ to $U > R$. We refer to such an event as a \textit{breakthrough}; a breakthrough increases the value of the investment project to
\begin{eqnarray} \label{at}
V_U(x)\equiv \sup_{\tau } \, \mathbf E_x \hskip 0.3mm [\mathrm e^{-r\tau } U (X_\tau )].
\end{eqnarray}
An immediate implication of (\ref{bof})--(\ref{at}) is that $V_U(x) > V_R(x)$, so that there is a surplus to be shared between the DM and the developers.

Developing the new technology involves a sunk cost $Z \in \mathbb R_+$, the same for each developer. We assume that $Z$ is the developers' private information; from the DM's perspective, $Z$ is drawn at time 0 from a known distribution, independently of $X$. Like the DM, each developer observes the evolution of $X$ and decides when to develop the new technology; she also observes the decisions of her competitor. To avoid coordination failures, we assume that, if both developers simultaneously attempt to develop the new technology, then only one of them, each with probability $1 \over 2$, can effectively have a breakthrough.\footnote{As in Katz and Shapiro (1987), this implicitly assumes that the developer losing the coin flip does not incur the cost $Z$. In Remark 1 below, we argue that dispensing with this public randomizing device would not significantly modify the DM's problem.}

The payoffs of the different parties can be described as follows. If the DM invests at time $\tau$ before a breakthrough occurs, then he obtains his stand-alone payoff $R(X_ \tau)$ and each developer obtains a zero payoff. If a breakthrough occurs before the DM invests in the project, at a time $\tau^B \leq \tau$ at which $X_{\tau^B} =x$, then the DM and the developer who has a breakthrough share the surplus $V_U(x) -V_R(x)$ according to the Nash (1950) bargaining solution. By himself, the DM can obtain $V_R(x)$ and, by herself, the successful developer can obtain zero. We deduce that their continuation payoffs at time $\tau^B$ are
\begin{eqnarray} \label{NASH}
G(x)  \equiv  \frac{1}{2} \,[V_{U} (x) + V_{R} (x)] \quad \mbox{and} \quad P(x)  \equiv  \frac{1}{2} \,[V_{U} (x) - V_{R} (x)],
\end{eqnarray}
respectively, while the other developer obtains a zero payoff.

For every strategy $\tau$ of the DM, these assumptions on payoffs give to the interaction between the developers the structure of a pure preemption game. Following Dutta and Rustichini (1993), we obtain that, in any pure-strategy subgame-perfect equilibrium, a breakthrough occurs at the first time at which the developers break even; that is, assuming that $P$ is strictly increasing and maps the state space for $X$ onto $\mathbb R_+$, we have
\begin{eqnarray*}
\tau^B = \tau_{X \geq Y} \equiv \inf \hskip 0.5mm \{ t\geq 0 : X_t \geq Y\},
\end{eqnarray*}
where $Y \equiv P^{-1}(Z)$. The upshot from this discussion is that, when $\tau \geq \tau_{X \geq Y}$, the DM pays $P(X_ { \tau_{X \geq Y}}) $ to a developer at time $\tau_{X \geq Y}$, which gives him the opportunity to benefit from a superior payoff $U(X_\tau)$ at time $\tau$. Thus the DM solves the following \pagebreak optimal stopping problem:
\begin{eqnarray} \label{11-2}
\sup _\tau \, { \mathbf E}_x \! \left[1_{\{\tau<\tau_{X \geq Y} \! \}}\, \mathrm e^{-r \tau} R(X_{\tau})+1_{\{\tau \geq \tau_{X \geq Y} \! \}}\! \left[ \mathrm e^{-r \tau}U(X_{\tau} )  - \mathrm e^{-r \tau_{X \geq Y}} P(X_ { \tau_{X \geq Y}}) \right] \right] \hskip -1mm.
\end{eqnarray}
Because $V_U(X_{\tau_{X\geq Y}}) - P(X_{\tau_{X\geq Y}}) = G(X_{\tau_{X\geq Y}})$ by (\ref{NASH}), the dynamic programming principle implies that the DM's problem (\ref{11-2}) can be rewritten as
\begin{eqnarray} \label{mainpb-intro}
\sup_{\tau} \,{ \mathbf E}_x \! \left[1_{\{\tau<\tau_{X \geq Y} \! \}}\, \mathrm e^{-r \tau} R(X_{\tau})+1_{\{\tau \geq \tau_{X \geq Y} \! \}}\, \mathrm e^{-r \tau_{X \geq Y}}G(X_{\tau_{X\geq Y}} )\right] \hskip -1mm.
\end{eqnarray}
This is the general form of the problem we shall deal with in this paper.

\paragraph{Remark}

Whenever a public randomization device is not available, the game can easily be modified to allow for simultaneous developments. In that case, it is convenient to assume that the DM and the developers share the surplus $V_U(x) -V_R(x)$ according to the Shapley (1953) value, which generalizes the Nash (1950) bargaining solution in the transferable-utility case. We distinguish two cases. Suppose first that only one developer, say, D1, develops the new technology at time $\tau^B$. Then D2 is effectively a null player, and hence she obtains a zero payoff; the DM's and D1's continuation payoffs at time $t$ are then given by (\ref{NASH}) as before. Suppose next that D1 and D2 simultaneously develop the new technology at time $\tau^B$. By himself, the DM can obtain $V_R(x)$ and, by themselves, D1 and D2 can obtain zero. Moreover, each developer does not increase the value of the coalition consisting of the DM and the other developer. We deduce that the continuation payoffs for the DM and each developer at time $\tau^B$ are
\begin{eqnarray} \label{SHAPLEY}
\overline G(x)  \equiv  \frac{1}{3} \,[2V_{U} (x) + V_{R} (x)] \quad \mbox{and} \quad \underline P(x)  \equiv  \frac{1}{6} \,[V_{U} (x) - V_{R} (x)],
\end{eqnarray}
respectively. Notice from \eqref{NASH} and \eqref{SHAPLEY} that $\overline G(x) > G(x)$ and $P(x) > \underline P(x)$, reflecting that the DM has a higher bargaining power when both developers have a breakthrough. The developers' strategies can be described using Riedel and Steg's (2017) adaptation of Fudenberg and Tirole's (1985) concept of extended mixed strategies to stochastic timing games. The equilibrium outcome is as follows. First, if $X_0 \leq P^{-1}(Z)$, then a single breakthrough occurs at time $\tau_{X \geq Y}$. Next, if $P^{-1} (Z) < X_0 < \underline P^{-1}(Z)$, then one or two simultaneous breakthroughs occur at time 0, depending on the realizations of the developers' extended mixed strategies. Finally, if $X_0 \geq \underline P^{-1}(Z)$, then two simultaneous breakthroughs occur at time 0. Thus the possibility of benefiting from two simultaneous breakthroughs can only increase the DM's payoff at time 0, which does not modify the solution to problem (\ref{mainpb-intro}).

\vskip 3mm

More generally, problem (\ref{mainpb-intro}) naturally arises in any situation in which a DM decides at each instant of time whether to settle for a basic payoff function $R$ or to wait until a continuation value function $G$ guaranteeing him a higher value than the option value associated to $R$ becomes available. As we will see in Section \ref{AMaFor}, because the time $\tau_{X \geq Y}$ at which this occurs is, by construction, the hitting time by $X$ of an unknown threshold $Y$, the Markovian formulation of problem (\ref{mainpb-intro}) leads to a two-dimensional stopping problem whose state variables are the process $X$ and its running maximum.

\section{Related Literature} \label{RelLit}

This paper is closely related to the literature on technology adoption under technological uncertainty pioneered by Balcer and Lippman (1984) and further developed by Weiss (1994), Farzin, Huisman, and Kort (1998), and Doraszelski (2004). We share with these authors the basic premise that the DM faces uncertainty about the arrival of a new technology or of an improvement on the current best practice. This implies that he has an option value of waiting before adopting the current state-of-the-art technology; for instance, Balcer and Lippman (1984) show that it is optimal for a firm to do so only if the technology it has in place lags behind by more than a certain amount. However, a key assumption of these models is that the value of the current state-of-the-art technology is fixed and known to the firm; thus, for instance, there is no output-price uncertainty. By contrast, we assume that this value---as well as the value of the new technology, should it become available---fluctuates randomly as a function of market conditions. Thus, in addition to technological uncertainty, the DM in our model faces cash-flow uncertainty, as in the standard real-options models of MacDonald and Siegel (1986) and Dixit and Pindyck (1994), or, in the context of technology adoption, as in the pure uncertain-profitability models of Jensen (1982), McCardle (1985), and Bhattacharya, Chatterjee, and Samuelson (1986). A contribution of our paper is to bring together these two classes of models of technology adoption.

In this respect, it is interesting to contrast our general investment problem and the investment-timing example of Section \ref{AMotExa} with recent developments of the real-options literature on investment under technological and cash-flow uncertainty. Papers in that literature, such as Alvarez and Stenbacka (2001), Huisman and Kort (2004), Murto (2007), Chronopoulos and Siddiqui (2015), and Chronopoulos and Lumbreras (2017), assume that technological breakthroughs are exogenous and arise independently of the realizations of the cash-flow process; a common assumption is that breakthroughs occur in a memoryless way, according to an independent Poisson process. This implies that current market conditions are the relevant Markov state variable for the firm's optimal investment policy. By contrast, breakthroughs are endogenous in our investment-timing example, and they arise in equilibrium when the market conditions become favorable enough to cover the developers' cost of introducing the new technology. This generates a rich two-dimensional dynamics involving current market conditions as well as their historic maximum, and leads to the stark prediction that investment in the stand-alone technology only takes place when market conditions deteriorate enough after having reached a maximum.

Grenadier and Weiss (1997) is closer to the present setup. However, they do not interpret the underlying stochastic process as a cash-flow process, but rather as describing the evolution of the state of technological progress; they assume that a breakthrough occurs when this process reaches a known threshold. By contrast, the threshold at which a breakthrough occurs in our investment-timing example is unknown to the DM as he does not observe the cost of introducing the new technology; as time goes by and the cash-flow process reaches new maximum values without a breakthrough occurring, the DM learns about this cost, becoming more pessimistic that he will eventually benefit from a breakthrough. Another difference is that the payoffs upon investing are functions of the current cash-flow in our model, whereas they are random variables independent of the state of technological progress in Grenadier and Weiss (1997).

From a technical viewpoint, this paper is related to the literature on two-dimensional stopping problems involving the running maximum of a
one-dimensional diffusion. Following the seminal contributions of Shepp and Shiryaev (1993), Dubins, Shepp, and Shiryaev (1994), Graversen and Peskir (1998), and Peskir (1998), an abundant literature in mathematical finance has used such models for the pricing of exotic options; see, for instance, Pedersen (2000), Guo and Shepp (2001), Dai and Kwok (2006), Guo and Zervos (2010), Ott (2014), and Rodosthenous and Zervos (2017). Compared to these studies, we do not a priori postulate an objective functional depending on the running maximum, as, for instance, in the analysis of Russian options; rather, we derive it from a natural investment problem using an appropriate change-of-variables formula. Moreover, a key distinctive feature of the Markovian formulation of our problem is that the objective functional does not involve the running maximum itself, but rather the integral of the discounted payoff from investing in the superior technology with respect to the increments of the maximum process, reflecting that the DM does not know at which threshold value of the underlying diffusion process a breakthrough will occur. Finally, most studies that provide an explicit solution to an optimal stopping problem involving the running maximum assume a specific functional form---typically, a geometric Brownian motion---for the underlying diffusion process; by contrast, our results hold for a rich class of time-homogenous diffusion processes.

\section{The Model} \label{TheMod}

In this section, we precisely formulate our general stopping problem and relate it to the model of investment under technological breakthroughs informally presented in Section \ref{AMotExa}.

\subsection{A General Stopping Problem}

Let $X \equiv (X_t)_{t \geq 0}$ be a one-dimensional time-homogeneous diffusion process defined over the canonical space $(\Omega,\CF,\mathbf P_x)$ of continuous trajectories with $X_0=x$ under $\mathbf P_x$, which is solution in law to the stochastic differential equation (SDE)
\begin{eqnarray} \label{eq1}
\mathrm dX_t = \mu(X_t) \, \mathrm dt  \, + \, \sigma(X_t) \, \mathrm dW_t, \quad t \geq 0,
\end{eqnarray}
driven by some Brownian motion $W \equiv (W_t)_{t \geq 0}$. The state space for $X$ is an interval $\mathcal I \equiv (\alpha, \beta)$, with $- \infty \leq \alpha < \beta \leq \infty$, and $\mu$ and $\sigma$ are continuous functions, with $\sigma>0$ over $\mathcal I$. We assume that $\alpha$ and $\beta$ are inaccessible (natural) endpoints for the diffusion. Therefore, $X$ is regular over $\mathcal I$ and the SDE \eqref{eq1} admits a weak solution that is unique in law. We also consider a random variable $Y$ with law $\mathbf Q$ taking values in $\mathcal I$ and independent of $X$. Overall, the relevant probability space for our analysis is the canonical product space $(\overline \Omega, \overline \CF, \overline {\mathbf P}_x) \equiv (\Omega \times \mathcal I , \CF \otimes \CB(\mathcal I), \mathbf P_x \otimes \mathbf Q)$, where $\CB(\mathcal I)$ is the Borel $\sigma$-field over $\mathcal I$. We denote by $\mathbf E_x$ and $\overline {\mathbf E}_x$ the expectation operators associated to $\mathbf P_x$ and $\overline {\mathbf P}_x$, respectively.

The DM observes the evolution of $X$, which he must stop at an appropriate time. Thus his strategy space is the set $\CT_X$ of all stopping times of the right-continuous filtration generated by $X$ over the canonical space; notice that the elements of $\CT_X$ can be identified to functions defined over $\overline \Omega$ and taking values in $\mathbb R_+ \cup \{ \infty\}$. Letting
\begin{eqnarray*}
\tau_{X \geq Y} \equiv \inf \hskip 0.5mm \{ t\geq 0 : X_t \geq Y\},
\end{eqnarray*}
we define the value for the DM of stopping $X$ at $\tau \in \CT_X$ as
\begin{eqnarray}
\overline J(x, \tau) \equiv  \overline{ \mathbf E}_x \! \left[1_{\{\tau<\tau_{X \geq Y} \! \}}\, \mathrm e^{-r \tau} R(X_{\tau})+1_{\{\tau \geq \tau_{X \geq Y} \! \}}\, \mathrm e^{-r \tau_{X \geq Y}}G(X_{\tau_{X\geq Y}} ) \right] \label{Jbar}
\end{eqnarray}
for some Borel functions $R$ and $G$ defined over $\mathcal I$.\footnote{\label{footnoteconvention}By convention, we let $f(X_{\tau}) \equiv 0$ over $\{ \tau = \infty \}$ for any Borel function $f$ and any random time $\tau$.} The interpretation of \eqref{Jbar} is that the DM chooses $\tau$ without knowing the realization of $Y$. If he stops the process $X$ before it reaches $Y$, then he is rewarded according to $R$, at time $\tau$. Otherwise, the game is stopped when $X$ reaches $Y$, and he is rewarded according to $G$, at time $\tau_{X \geq Y}$. Our objective is to solve the following optimal stopping problem:
\begin{eqnarray} \label{mainpb}
\overline V (x) \equiv  \sup_{\tau \in \CT_X} \overline J(x, \tau).
\end{eqnarray}
Before we detail the technical assumptions we impose on the primitives of the model, it is helpful to relate problem (\ref{mainpb}) to the investment problem with technological breakthroughs (\ref{11-2}), which we used to motivate our analysis in Section \ref{AMotExa}. In problem (\ref{11-2}), the stopping time $\tau$ belongs to $\CT_{X, X \geq Y}$, the set of all stopping times of the filtration $(\CG_t)_{t\geq 0}$ over $\overline \Omega$ defined by $\CG_t \equiv \sigma( X_s,  1_{\{\tau_{X \geq Y} \leq s\}}; s \leq t)$ for all $t \geq 0$. That is, the information that accrues to the DM up to any time $t$ is the evolution of $X$ up to time $t$, as well as the breakthrough time at which $X$ reaches $Y$, should this happen before time $t$. For the sake of completeness, we verify in Appendix B that the dynamic programming principle applies to problem (\ref{11-2}), which allows us to rewrite it under the general form (\ref{mainpb}).

\subsection{Technical Assumptions} \label{TecAss}

We first recall useful properties of the solution $X$ to the SDE (\ref{eq1}). We next detail the assumptions on the payoff functions $R$ and $G$ and on the law $\mathbf Q$ of the random variable $Y$ under which we solve problem (\ref{mainpb}). We also emphasize useful properties of the following auxiliary optimal stopping problem:
\begin{eqnarray} \label{sa}
V_R(x) \equiv \sup_{\tau \in \CT_X} \mathbf E_x \hskip 0.3mm[ \mathrm e^{-r\tau } R(X_\tau)],
\end{eqnarray}
which plays an important role in our analysis. Intuitively, (\ref{sa}) corresponds to the stand-alone investment problem in which the DM cannot benefit from a technological breakthrough.

\paragraph{Properties of the Diffusion $X$}

The infinitesimal generator of the diffusion $X$ is defined for functions $u \in \mathcal C^2(\mathcal I)$ by
\begin{eqnarray}
\CL u (x) \equiv \mu(x) u'(x) + \frac{1}{2}\,\sigma^2(x) u''(x), \quad x \in \mathcal I. \label{defmathcalL}
\end{eqnarray}
That $\sigma >0$ over $\mathcal I$ ensures that the equation $\CL u - ru=0$ admits a two-dimensional space of solutions in $\mathcal C^2(\mathcal I)$, spanned by two positive fundamental solutions $h_1$ and $h_2$, respectively strictly increasing and strictly decreasing, that are uniquely defined up to a linear transformation. By Abel's theorem, the ratio
\begin{eqnarray}
\gamma \equiv  \frac{h'_1(x)h_2(x) -h_1(x) h_2'(x) }{S'(x)} >0 \label{wronskian}
\end{eqnarray}
of the Wronskian of $h_1$ and $h_2$ and of the derivative of the scale function of the diffusion $X$, which is uniquely defined up to an affine transformation by
\begin{eqnarray}
S(x)\equiv \int_c^x \exp\left(-\int_c^y \frac{2 \mu(z)}{\sigma^2(z)}\, \mathrm dz \right)\mathrm dy, \quad x \in \mathcal I \label{scale}
\end{eqnarray}
for some fixed $c \in \mathcal I$, is a constant independent of $x$. Because the boundaries $\alpha $ and $\beta$ of $\mathcal I$ are natural, we know in particular that
\begin{eqnarray}
\lim_{x \to \alpha^+} h_1(x)=0, \quad \lim_{x \to \beta^-}h_1(x)=\infty, \quad \lim_{x \to \alpha^+} h_2(x)=\infty, \quad \lim_{x \to \beta^-}h_2(x)=0. \label{hboundaries}
\end{eqnarray}
Furthermore, defining the hitting time $\tau(y) \equiv \inf \hskip 0.5mm \{ t \geq 0 : X_t = y \}$ for all $y \in \mathcal I$, the mapping
\begin{eqnarray}
\label{laplace} (x,y) \mapsto \mathbf E_x \hskip 0.3mm [\mathrm e^{-r\tau(y)}] = \left\{ \begin {matrix} \frac{h_1(x)}{h_1(y)} & \text{if} & x \leq y, \\  \frac{h_2(x)}{h_2(y)} & \text{if} & x > y, \end{matrix} \right.
\end{eqnarray}
is continuous and $\mathcal C^2$ over $\{(x,y) \in\mathcal I \times \mathcal I: x \not = y\}$.

\paragraph{Assumptions on the Payoff Functions $R$ and $G$}

We assume that $R \in \mathcal C^2(\mathcal I)$, and that it satisfies
\begin{itemize}

\item[\bf A1]

For each $x \in \mathcal I$, $\mathbf E _x \hskip 0.3mm  [ \sup_{t \geq 0} \mathrm e^{-rt} | R(X_t) | ] < \infty$.

\item[\bf A2]

For each $x \in \mathcal I$, $\lim_{t \to \infty} \mathrm e^{-rt} R(X_t) = 0$, $\mathbf P_x$-almost surely.

\item[\bf A3]

There exists $x_0 \in \mathcal I$ such that ${\cal L}R - rR >0$ over $(\alpha, x_0)$ and ${\cal L}R - rR < 0$ over $( x_0, \beta)$.

\end{itemize}
A1 guarantees that the family $(\mathrm e^{-r  \tau} R(X_{ \tau}))_{\tau \in {\cal T}_X}$  is uniformly integrable. A1--A2 imply the useful growth property
\begin{eqnarray} \label{gp}
\lim_{x \to \alpha^+} \frac{R(x)}{h_2(x)} = \lim_{x \to \beta^-} \frac{R(x)}{h_1(x)} = 0
\end{eqnarray}
and are in line with the convention made in Footnote \ref{footnoteconvention}. A3 guarantees that the optimal stopping region $\{x \in \mathcal I: V_R(x) = R(x)\}$ for the stand-alone optimal stopping problem (\ref{sa}) is of the form $[x_R,\beta)$ for some threshold $x_R > x_0$, so that
\begin{eqnarray}
V_R(x) = \left\{ \begin{array}{lll} \frac{h_1(x)}{h_1(x_R)} \, R(x_R )& \text{if} & x < x_R, \\  R(x) & \text{if} & x \geq x_R, \end{array} \right. \label{VR}
\end{eqnarray}
and the smooth-fit property applies at $x_R$, that is, $R'(x_R)=\frac{h'_1(x_R)}{h_1(x_R)}\,R(x_R)$ (Peskir and Shiryaev (2006), Dayanik and Karatzas (2003, Corollary 7.1)). It follows from standard optimal stopping theory that $(\mathrm e^{-rt} V_R(X_t))_{t \geq 0}$ is a supermartingale and that ${\cal L}V_R - rV_R \leq 0 $ over $\mathcal I \setminus \{x_R\}$. The following lemma holds.

\begin{lem} \label{lem_sign}
$V_R  >0$ over $\mathcal I$ and $R >0$ over $[x_R, \beta).$
\end{lem}

We assume that $G \in \mathcal C^1( \mathcal I)$, that $G$ is piecewise $\mathcal C^2$ over $\mathcal I$, and that it satisfies
\begin{itemize}

\item[\bf A4]

For each $x \in \mathcal I$, $\mathbf E_x \hskip 0.3mm [ \sup_{t \geq 0} \mathrm e^{-rt} G(X_t)] < \infty $.

\item[\bf A5]

For each $x \in \mathcal I$, $\lim_{t \to \infty} \mathrm e^{-rt} G(X_t) = 0$, $\mathbf P_x$-almost surely.

\item[\bf A6]

$G > V_R$ over $\mathcal I$.

\item[\bf A7]

${\cal L}G - rG \leq 0  $ everywhere $G''$ is defined.

\end{itemize}
From (\ref{sa}), A6, and Lemma \ref{lem_sign}, we have $G > R\vee 0$ over $\mathcal I$; hence A4 guarantees that the family $(\mathrm e^{-r  \tau} G(X_{ \tau}))_{\tau \in {\cal T}_X}$ is uniformly integrable. A4--A5 imply the useful growth property
\begin{eqnarray} \label{gp'}
\lim_{x \to \alpha^+} \frac{G(x)}{h_2(x)} = \lim_{x \to \beta^-} \frac{G(x)}{h_1(x)} = 0.
\end{eqnarray}
The interpretation of A6--A7 is that $G$ dominates the value function $V_R$ of the stand-alone optimal stopping problem (\ref{sa}) and that it incorporates itself the solution to an optimal stopping problem, so that $(\mathrm e^{-rt} G(X_t))_{t \geq 0}$ is a supermartingale. In the investment problem with technological breakthroughs (\ref{11-2}), this occurs because, at time $\tau_{X \geq Y}$, the developers allow the DM to substitute to the stand-alone technology, with payoff function $R$, a more efficient technology, with payoff function $U > R$. The function $G$ results from a subsequent optimal stopping problem involving the superior payoff function $U$, as shown by (\ref{at})--(\ref{NASH}).

\paragraph{Assumptions on the Distribution of $Y$}

Recall that the random variable $Y$ takes values in $\mathcal I$ and is independent of $X$. We further assume that its law $\mathbf Q$ satisfies
\begin{itemize}

\item[\bf A8]

$\mathbf Q$ has locally Lipschitz density $f >0$ over $\mathcal I$ with respect to Lebesgue measure.

\end{itemize}
We denote by $F$ the cumulative distribution function of $Y$.

\vskip 3mm

In line with the real-options literature (Dixit and Pindyck (1994)), a natural specification of the model consists in letting $X$ follow a geometric Brownian motion with drift $\mu<r$ and volatility $\sigma >0$, and in letting the payoff functions be given by $R(x) \equiv x-I$ and $U(x) \equiv \kappa x -I$; here $I$ is a positive investment cost that the DM must incur to obtain a cash-flow $X$ (under the stand-alone technology) or $\kappa X$ (under the more efficient technology), where $\kappa> 1$. Finally, we assume that the development cost $Z$ is drawn from a distribution with locally Lipschitz density $f_Z >0$ over $(0,\infty)$ with respect to Lebesgue measure. We verify in Appendix C that this specification satisfies A1--A8.

\section{A Markovian Formulation} \label{AMaFor}

We now turn to the analysis of our main problem (\ref{mainpb}), for which we first give a convenient Markovian formulation. Given $m \in [x, \beta)$, define the maximum process $M \equiv (M_t)_{t \geq 0}$ by
\begin{eqnarray*}
M_t= m \vee \sup_{s \leq t} \hskip 0.3mm X_s
\end{eqnarray*}
for all $t \geq 0$, so that the pair $(X,M)$ defines a continuous Markov process starting at $(x,m)$. We denote by $\mathbf P_{x,m}$ the law of this Markov process over $(\Omega, \mathcal F) \otimes (\Omega, \mathcal F)$ and by $\mathbf E_{x,m}$ the corresponding expectation operator. Observe that any stopping time in the set $\CT_{X,M}$ of all stopping times of the right-continuous filtration generated by $(X,M)$ over $(\Omega, \mathcal F)\otimes (\Omega, \mathcal F)$ is $\mathbf{P}_{x,m}$-almost surely equal to a stopping time in $\CT_X$.

\begin{pro} \label{markovian}
The function
\begin{eqnarray} \label{mainproblem}
V(x,m)\equiv \sup_{\tau \in \CT_{X,M}} \mathbf E _{x,m} \! \left[[1-F(M_{\tau})] \, \mathrm e^{-r \tau}R(X_{\tau}) +\int_{0}^\tau \mathrm e^{-r t}G(M_t)f(M_t) \, \mathrm dM_t \right]
\end{eqnarray}
is well-defined and
\begin{eqnarray}
\overline V(x) = V(x,x) + F(x) G(x) \label{overlineVVFG}
\end{eqnarray}
for all $x \in \mathcal I$.
\end{pro}

To grasp the intuition for this result, consider a static version of our problem in which the DM has no decision to take and all uncertainty is resolved immediately at time $0$. Then the DM obtains $G(x)$ with probability $F(x)$, as in (\ref{overlineVVFG}), and $R(x)$ with probability $1 - F(x)$, as in (\ref{mainproblem}) for $\tau = 0$. The integral with respect to the maximum process in (\ref{mainproblem}) should thus be interpreted as the added value of postponing investment in the hope of a technological breakthrough. Specifically, because breakthroughs only occur when the process $X$ reaches new maximum values, the probability that a breakthrough occurs during the time interval $(0,\tau)$ is $F(M_\tau)-F(m)$, and the expected discounted value of such breakthroughs to the DM, compounded over the increments of $M$ over $(0,\tau)$, is given by
\begin{eqnarray*}
{\int_{0}^\tau \mathrm e^{-r t}G(M_t)f(M_t) \, \mathrm dM_t \over F(M_\tau)-F(m)}.
\end{eqnarray*}
By contrast, the downside risk of postponing investment is that, with probability $1- F(M_\tau )$, no breakthrough may occur before time $\tau$, in which case the DM may well end up with a payoff $R(X_ \tau) < R(x)$.

\bigskip

\noindent \textbf{Proof.} Because the diffusion $X$ is regular over $\mathcal I$, we know that, for each $y \in [x, \beta)$, the hitting time $\tau(y) = \inf \hskip 0.5mm \{ t\geq 0 : X_t = y \}$ is finite with positive probability under $\mathbf P_x$. Notice that $\tau_{X \geq Y}=\tau(Y)$ over $\{Y >x\}$ and $\tau_{X \geq Y} =0$ over $\{Y \leq x\}$. By Fubini's theorem, we have
\begin{align}
\overline J(x,\tau) &= \mathbf E_{x}\! \left[ \int_{x}^\beta \! \left[ 1_{\{\tau<\tau(y)\}} \, \mathrm e^{-r \tau}R(X_{\tau})+1_{\{\tau \geq \tau(y)\}} \, \mathrm  e^{-r \tau(y)}G(y)\right] \! \mathbf Q(\mathrm d y)+ F(x) G(x) \right] \notag \allowdisplaybreaks
\\
&= \mathbf E_{x,x}\! \left[ \int_{M_\tau}^{\beta} f(y)\, \mathrm dy \; \mathrm e^{-r \tau}R(X_{\tau})+\int_{x}^{M_\tau}\mathrm e^{-r \tau(y)}G(y ) f(y)  \, \mathrm dy  \right] + F(x) G(x)\notag
\\
&= \mathbf E_{x,x} \!\left[ [1-F(M_{\tau})] \, \mathrm e^{-r \tau}R(X_{\tau})+\int_{x}^{M_\tau}\mathrm e^{-r \tau(y)}G(y) f(y) \, \mathrm dy \right]  + F(x) G(x). \label{trucmalin}
\end{align}
By A1, we have
\begin{eqnarray*}
\mathbf E_{x,m}\! \left[[1-F(M_{\tau})]\, \mathrm e^{-r \tau}R(X_{\tau}) \right] < \infty
\end{eqnarray*}
for all $\tau \in {\cal T}_{X,M}$. The remainder of the proof of Proposition \ref{markovian} relies on the following lemma, which establishes a change-of-variables formula that clarifies the dependence of the payoff with respect to the maximum process $M$.

\begin{lem}
\label{lemchange-of-time} For each $\tau \in {\cal T}_{X,M},$
\begin{eqnarray}\label{change-of-time}
\int_{m}^{M_\tau} \mathrm e^{-r \tau(y)}G(y) f(y) \, \mathrm dy= \int_{0}^\tau \mathrm e^{-r t}G(M_t)f(M_t) \, \mathrm dM_t
\end{eqnarray}
$\mathbf P_{x,m}$-almost surely.
\end{lem}

\noindent \textbf{Proof.} Fix some $\omega \in \Omega$, and consider the continuous nondecreasing mapping $t \mapsto M_t$ and its right-continuous inverse
\begin{eqnarray*}
C_y \equiv \inf \hskip 0.3mm \{ t\geq 0 :  M_t > y\}, \quad y \in[m,\beta),
\end{eqnarray*}
with $\inf \hskip 0.5mm \emptyset  = \infty$ by convention. By construction, $C_{y-}=\tau(y)$ for all $y \in (m, \beta)$. Consider the Borel function
\begin{eqnarray*}
g(t) \equiv 1_{\{0<t \leq \tau\}} \, \mathrm e^{-rt}G(M_t)f(M_t), \quad t \geq 0,
\end{eqnarray*}
with $g(\infty) \equiv 0$ by convention, which is nonnegative as $G >0$ over $\mathcal I$. According to the change-of-variables formula for Stieltjes integrals (Revuz and Yor (1999, Chapter 0, Proposition 4.9)), we have
\begin{eqnarray}\label{eq_change_time}
\int_{[0,\infty)} g(t) \, \mathrm d M_t = \int_m^\beta g(C_y) \, \mathrm  dy
\end{eqnarray}
whenever these integrals are well-defined. To check that this is the case, notice first that, because $M$ is continuous, $y= M_{\tau(y)}$ for all $y \in [m, \beta)$; moreover, because $C$ is right- continuous and nondecreasing, $\tau(y) =C_{y-}= C_{y}$ for all $y \in [m, \beta)$ outside of a countable set. It follows that
\begin{align}
\int_m^{\beta} g(C_y) \, \mathrm dy  &= \int_m^\beta 1_{\{0< C_y \leq \tau\}} \, \mathrm e^{-r C_y}G(M_{C_y})f(M_{C_y})  \, \mathrm dy \notag
\\
&= \int_{m}^{M_\tau} \mathrm e^{-r \tau(y)}G(y) f(y)   \, \mathrm dy.  \label{271}
\end{align}
We claim that the quantity (\ref{271}) is $\mathbf P_{x,m}$-almost surely finite. Indeed, taking expectations and using Fubini's theorem, we have
\begin{align*}
\mathbf E_{x,m} \!  \left[ \int_{m}^{M_\tau} \mathrm e^{-r \tau(y)}G(y) f(y)  \, \mathrm dy \right] & \leq   \int_m^\beta \mathbf E_x  \! \left[\mathrm e^{-r \tau(y)}\right ] \! G(y) f(y) \, \mathrm dy \allowdisplaybreaks
\\
& \leq [1 - F(m)] \sup_{\tau\in \mathcal T_X} \mathbf E_x \! \left[ \mathrm e^{-r \tau} G(X_\tau)  \right]\hskip -1mm,
\end{align*}
which is finite by A4. Therefore, the integrals in \eqref{eq_change_time} are well-defined, as claimed. To conclude the proof of Lemma \ref{lemchange-of-time}, simply observe that the left-hand side of \eqref{eq_change_time} is equal to the right-hand side of (\ref{change-of-time}) and that, in line with \eqref{271}, the right-hand side of \eqref{eq_change_time} is equal to the left-hand side of (\ref{change-of-time}). The result follows. \hfill $\blacksquare$

\bigskip

We are now ready to complete the proof of Proposition \ref{markovian}. By Lemma \ref{lemchange-of-time}, the value function $V$ of problem \eqref{mainproblem} is well-defined and, by (\ref{trucmalin})--(\ref{change-of-time}), we have
\begin{align*}
\overline{V}(x) &= \sup_{\tau \in \mathcal T_X} \overline J(x,\tau)
\\
&= \sup_{\tau \in \mathcal T_{X,M}} \mathbf E_{x,x} \! \left[ [1-F(M_{\tau})] \, \mathrm e^{-r \tau}R(X_{\tau})+\int_{x}^{M_\tau} \mathrm e^{-r \tau(y)}G(y) f(y)  \, \mathrm dy \right]  + F(x) G(x) \allowdisplaybreaks
\\
&= \sup_{\tau \in \mathcal T_{X,M}} \mathbf E_{x,x} \! \left[ [1-F(M_{\tau})] \, \mathrm e^{-r \tau}R(X_{\tau})+ \int_{0}^\tau \mathrm e^{-r t}G(M_t)f(M_t) \, \mathrm dM_t \right]  + F(x) G(x),
\end{align*}
which is \eqref{overlineVVFG} by \eqref{mainproblem}. Hence the result. \hfill $\blacksquare$

\bigskip

Our initial problem (\ref{mainpb}) hence reduces to an optimal stopping problem for the two- dimensional Markov process $(X,M)$ over the state space $\mathcal J \equiv \{(x,m) \in \mathcal I \times \mathcal I : m \geq x\}$. Compared to similar problems involving the maximum process so far considered in the literature, a distinctive feature of problem (\ref{mainproblem}) is that the DM's payoff function features an integral with respect to the maximum process.

\section{The Main Theorem} \label{TheMaiThe}

In this section, we first heuristically derive a variational system for the value function $V$ of problem (\ref{mainproblem}). Our main theorem then states that this system has a unique solution, which coincides with $V$, and expresses the optimal stopping time for (\ref{mainproblem}) in terms of a free boundary in the domain $\mathcal J$.

\subsection{A Heuristic Derivation}

\paragraph{The Stopping Region}

We start with three educated guesses about the optimal stopping region $\mathcal S \subset \mathcal J$ for problem (\ref{mainproblem}), which are in line  with Peskir's (1998) classic analysis of the stopping problem for the maximum process.

First, because $G > V_R$ by A6, a lower bound for $V$ is obtained by stopping $X$ at the optimal threshold $x_R$ for the stand-alone optimal stopping problem (\ref{sa}). Intuitively, it is thus suboptimal for the DM to stop $(X,M)$ before $X$ reaches $x_R$. That is, we guess
\begin{eqnarray*}
\mathcal S \subset \{(x,m) \in \mathcal J : x \geq x_R\}.
\end{eqnarray*}
Second, the maximum process $M$ can increase only when $(X, M)$ hits the diagonal $\mathcal D \equiv \{(x,m ) \in \mathcal J : x =m\}$ of $\mathcal J$, at which point the DM's payoff may jump upwards. Intuitively, it is thus suboptimal for him to stop $(X, M)$ over $\mathcal D$. That is, we guess
\begin{eqnarray*}
\mathcal S \cap \mathcal D = \emptyset.
\end{eqnarray*}
Third, after $(X, M)$ hits $\mathcal D$, it is costly for the DM to let it run horizontally too far to the left of $\mathcal D$ because of the time needed to reach a new value of $M$. Intuitively, because of discounting, the opportunity cost of delaying action increases when $(X,M)$ moves away from $\mathcal D$. That is, we guess
\begin{eqnarray*}
\mathcal S= \{ (x,m) \in \mathcal J : \underline m \leq m <\beta  \text{ and } x_R \leq x \leq  b(m)\}
\end{eqnarray*}
for some $\underline m \in ( x_R, \beta)$ and some function $b:[\underline m, \beta) \to [x_R, \beta)$ that satisfies $b(\underline m)= x_R$, $b(m) < m $ for all $m \in [\underline m, \beta)$, and the limit condition
\begin{eqnarray*}
\lim_{m \to \beta^-}b(m) = \beta.
\end{eqnarray*}
The fact that the free boundary $x= b(m)$ is defined for $x \geq x_R$ and $m \geq \underline m$ follows naturally from our investment problem and is specific to our model. We will see that the limit condition turns out to be a consequence of the growth property (\ref{gp'}). This reflects the idea that the DM should be reluctant to let the process $(X,M)$ run too far from $\mathcal D$ as $X$ grows large.

\paragraph{Dynamic Programming}

The dynamic programming principle for (\ref{mainproblem}) relies on the key observation that $M$ does not increase as long as $(X,M)$ stays away from $\mathcal D$. Hence the infinitesimal generator of $(X, M)$ restricted to functions $u \in \mathcal C^{2,1}(\mathrm {int} \, \mathcal J)$ coincides over $\mathrm {int} \, \mathcal J$ with the infinitesimal generator of $X$,
\begin{eqnarray*}
\CL u (x,m) \equiv \mu(x) \, \frac{\partial u}{\partial x}\,(x,m)+ \frac{1}{2}\, \sigma^2(x) \, \frac{\partial^2 u}{\partial x^2} \, (x,m), \quad (x,m)  \in \mathrm {int} \, \mathcal J.
\end{eqnarray*}
Letting $\mathcal C \equiv \mathcal J \setminus \mathcal S$ be the continuation region, the dynamic programming principle then states that, provided $V$ is $\mathcal C^{2,1}$ over $\mathrm {int}\, \mathcal C \subset \mathrm {int} \, \mathcal J$, we have
\begin{eqnarray*}
\CL V(x,m)-rV(x,m) =0, \quad (x,m) \in \mathrm {int}\, \mathcal C.
\end{eqnarray*}
The dynamic programming principle is at this stage only an educated guess, which we will ultimately confirm by a verification argument.

\paragraph{Value-Matching}

As usual, the value-matching condition
\begin{eqnarray*}
V(x,m) = [1 - F(m)] R(x) , \quad (x,m) \in \mathcal S
\end{eqnarray*}
pins down the value function on the stopping region.

\paragraph{Smooth-Fit}

We conjecture that the value function satisfies the smooth-fit property along horizontal lines at the free boundary $x= b(m)$:
\begin{eqnarray*}
\label{sp710} {\partial V \over \partial x}\,(b(m),m) = [1-F(m)]R'(b(m)),\quad  m\in[\underline m,\beta).
\end{eqnarray*}
Like the dynamic programming principle, the smooth-fit property is at this stage only an educated guess, which we will ultimately have to verify using an appropriate candidate for the free boundary.

\paragraph{Neumann Condition}

The Neumann condition expresses the behavior of the value function at $\mathcal D$, where the process $(X,M)$ undergoes a normal reflection. In our setting, the Neumann condition takes the form
\begin{eqnarray*}
\frac{\partial V}{\partial m}\,(m,m) = -f(m)G(m),\quad  m \in \mathcal I,
\end{eqnarray*}
which can be heuristically derived by observing that, from \eqref{mainproblem} and \eqref{change-of-time},
\begin{eqnarray*}
V(m,m)\equiv \sup_{\tau \in \CT_{X,M}} \mathbf E _{m,m} \! \left[[1-F(M_{\tau})]\, \mathrm e^{-r \tau}R(X_{\tau}) +\int_{m}^{M_\tau} \mathrm e^{-r \tau(y)}G(y) f(y) \, \mathrm dy \right]\hskip -1mm .
\end{eqnarray*}
Intuitively, starting from a point $(m,m) \in \mathcal D$, a marginal increase in the second argument brings bad news to the DM. Indeed, noticing that $f(m)\, \mathrm dm$ corresponds to the probability that a breakthrough occurs in $[m, m + \mathrm dm]$, the term $f(m)G(m)$ on the right-hand side of the Neumann condition represents the expected foregone payoff for the DM if $X$ reaches a new maximum value yet no breakthrough occurs. It should be noted that this heuristic derivation incorporates our guess that the free boundary $x=b(m)$ does not cut $\mathcal D$.

\paragraph{Boundary Condition}

Finally, because the lower endpoint $\alpha$ of $\mathcal I$ is inaccessible for $X$ and $\lim_{x \to \alpha^+} V_R(x) = 0$ by \eqref{hboundaries} and \eqref{VR}, we conjecture that the value function vanishes at $\alpha$,
\begin{eqnarray*}
\lim_{x \to \alpha^+} V(x,m) = 0, \quad  m \in \mathcal I.
\end{eqnarray*}
Overall, we are led to find a threshold $\underline m \in (x_R, \beta)$, a function $b: [\underline m, \beta) \to [x_R,\beta)$, and a function $W : \mathcal J \to \mathbb R$ that jointly satisfy the following variational system (VS):
\begin{align}
b(\underline m) &=  x_R, \label{ic} \allowdisplaybreaks
\\
x_R \hskip 0.4mm \leq \hskip 0.4mm b(m) &< m, \quad m \in [\underline m, \beta), \label{lalala}
\\
\lim_{m \to \beta^-} b(m) &= \beta,  \label{tc}
\\
\CL W(x,m)-rW(x,m) &= 0, \quad (x,m) \in \mathrm {int}\, \mathcal C,  \label{vs1}
\\
W(x,m) &= [1-F(m)]R(x), \quad  (x,m) \in \mathcal S ,  \label{vs2}
\\
{\partial W \over \partial x}\,(b(m),m) &= [1-F(m)]R'(b(m)), \quad m\in[\underline m,\beta), \label{vs3}
\\
\frac{\partial W}{\partial m}\, (m,m) &= -f(m) G(m),\quad  m \in \mathcal I,  \label{vs4}
\\
\lim_{x \to \alpha^+} W(x,m)  &= 0, \quad  m \in \mathcal I,  \label{vs5}
\end{align}
where $\mathcal S \equiv \{ (x,m) \in \mathcal J : \underline m \leq m <\beta  \text{ and } x_R \leq x \leq  b(m)\}$ and $\mathcal C \equiv \mathcal J \setminus \mathcal S$.

\paragraph{Remark}

As in the literature on two-dimensional stopping problems involving the running maximum of a one-dimensional diffusion, an important feature of our variational system is that the smooth-fit condition is not sufficient to characterize the solution that coincides with the value function of the optimal stopping problem. To address this issue, the classical method of Shepp and Shiryaev (1993) and Dubins, Shepp, and Shiryaev (1994) consists in imposing a well-chosen growth condition on the free boundary or on the solution of the variational system; this is analogous to a transversality condition and allows them to pin down the optimal free boundary. The corresponding solution to the ODE turns out to be the largest solution which stays strictly below $\mathcal D$, a property referred to by Peskir (1998) as the \textit{maximality principle}. In our problem, it is not straightforward to guess the appropriate growth condition, nor to apply the maximality principle. Instead, the simple limit condition \eqref{tc} turns out to characterize the optimal free boundary among all candidate solutions.

\subsection{A Formal Statement}

Our central theorem can now be stated as follows.

\begin{thm} \label{maintheo}
(VS) admits a unique solution $(\underline m,b,W)$ in $(x_R, \beta) \times \mathcal C^1([\underline m, \beta)) \times \mathcal V,$ where
\begin{align*}
\mathcal V &\equiv \mathcal C^0(\mathcal J) \cap {\cal C}^1(\mathcal J \setminus \mathcal J_{\underline m}) \cap {\cal C}^{2,1} (\mathcal J_1)\cap {\cal C}^{2,1} (\mathcal J_2)\cap {\cal C}^{2,1} (\mathcal J_3)  \cap {\cal C}^{2,1} (\mathcal S) ,
\\
\mathcal J_{\underline m} &\equiv \{(x,\underline m) \in \mathcal J: x \in (\alpha, \underline m ]\},
\\
\mathcal J_1 &\equiv \{ (x,m) \in \mathcal J : m \in [\underline m, \beta) \mbox{ \rm{and} } x \leq x_R \},
\\
 \mathcal J_2 &\equiv \{ (x,m) \in \mathcal J : m \in [\underline m, \beta) \mbox{ \rm{and} } x \geq b(m)) \},
\\
\mathcal J_3 &\equiv \{ (x,m) \in \mathcal J : m \in (\alpha, \underline m] \}.
\end{align*}
Moreover$,$
\begin{itemize}

\item[(i)]

the function $b$ is strictly increasing$;$

\item[(ii)]

the function $W$ coincides with the value function $V$ of problem (\ref{mainproblem})$;$

\item[(iii)]

the stopping time
\begin{eqnarray*}
\tau_b \equiv \inf \hskip 0.5mm \{ t\geq 0  : M_t \geq \underline m \mbox{ \rm{and} } x_R \leq  X_t \leq b(M_t) \}
\end{eqnarray*}
is optimal for problem (\ref{mainproblem})$,$ that is$,$
\begin{eqnarray*}
V(x,m) = \mathbf E_{x,m} \! \left[[1-F(M_{\tau_b})]\, \mathrm e^{-r \tau_b } R(X_{\tau_b }) +\int_{0}^{\tau_b }\mathrm e^{-r t}G(M_t)f(M_t)\, \mathrm dM_t  \right]
\end{eqnarray*}
for all $(x,m) \in \mathcal J$.

\end{itemize}
\end{thm}

\paragraph{Remark}

Notice that a function in $\mathcal V$ need not be globally $C^1$ over $\mathcal J$. Indeed, it will turn out that the partial derivative with respect to $m$ of the function $W$ part of the solution to (VS) is not continuous at the points $(x, \underline m)$ for $x < \underline m$. This relates to the fact that the DM is willing to invest in the stand-alone technology only after $X$ reaches $\underline m$. We elaborate on this distinctive feature of our solution in Proposition \ref{propertyV1} of Section \ref{Discussion}.
\noindent \setcounter{figure}{0}
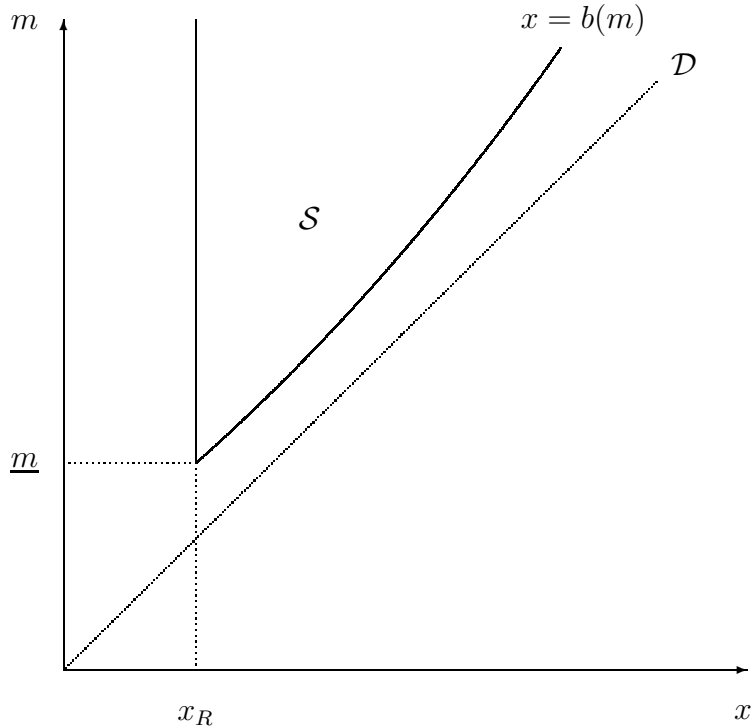
\begin{figure}
{\normalsize \noindent \unitlength=1.50mm \special{em:linewidth
0.4pt} \linethickness{0.4pt} \ \ \ \ \ \ \ \ \
\begin{picture}(140.00,67.5)
\put(18.4,10){\vector(1,0){60}}
\put(18.4,10){\vector(0,1){57.5}}
\put(78,6.28){\makebox(0,0)[cc]{$x$}}
\put(15,67.3){\makebox(0,0)[cc]{$m$}}
\bezier{242}(18.4,10)(44.4,36)(70.4,62)
\put(72.8,63.6){\makebox(0,0)[cc]{$\mathcal D$}}
 \put(40,50){\makebox(0,0)[cc]{$\mathcal S$}}


\bezier{1000}(30,28.29)(47,43.5)(62,65)
\bezier{1000}(30,28.29)(30,49.145)(30,67.5)
\put(64,67.5){\makebox(0,0)[cc]{$x= b(m)$}}
\put(15,28.29){\makebox(0,0)[cc]{$\underline m$}}
\put(30,6){\makebox(0,0)[cc]{$x_R$}}

\bezier{36}(30,28.29)(30,19.145)(30,10)
\bezier{23}(30,28.29)(24.34,28.29)(18.4,28.29)

\end{picture}}\vskip -7mm
\caption{The optimal stopping region.}
\end{figure}

\vskip 3mm

The optimal stopping region is illustrated on Figure 1. The proof of Theorem \ref{maintheo} consists of two parts.

The first part of the proof is purely analytical, and consists in showing that (VS) admits a solution $(\underline m, b, W)$ that satisfies the required regularity conditions. The function $b$ is characterized as the solution to an ordinary differential equation (ODE), and the function $W$ is explicit given $b$. This part of the proof is provided in Section \ref{analytical}.

The second part of the proof follows the theory of optimal stopping. We exploit two properties, which we show are satisfied by any function $W$ part of a solution to (VS). The first property is that $W$ is bounded below and above as follows:
\begin{eqnarray*}
[1-F(m)]R(x) \vee 0 \leq W(x,m) < [1-F(m)] G(x)
\end{eqnarray*}
for all $(x,m) \in \mathcal J$; this notably implies a useful uniform integrability property. The second property is that $W$ is superharmonic; that is, for any stopping time $\tau \in \mathcal T_{X,M}$,
\begin{eqnarray*}
W(x,m)\geq \mathbf E_{x,m} \! \left[\mathrm e^{-r \tau} W (X_\tau,M_\tau) +\int_{0}^\tau \mathrm e^{-r t}G(M_t)f(M_t) \, \mathrm dM_t \right]
\end{eqnarray*}
for all $(x,m) \in \mathcal J$. We then argue that these two properties together imply that $W$ must coincide with the value function $V$ of problem (\ref{mainproblem}) and that the stopping time $\tau_b$ is optimal for this problem. Incidentally, this also shows that (VS) admits a unique solution. This part of the proof is provided in Section \ref{section_verification}.

\section{Analysis of the Variational System} \label{analytical}

The central result of this section is that (VS) admits a solution.

\begin{pro} \label{main3}
There exist $\underline m \in (x_R,\beta),$ a strictly increasing function $b\in \mathcal C^1([\underline m, \beta)),$ and a function $W \in \mathcal V$ such that $(\underline m, b,W)$ is a solution to (VS).
\end{pro}

It should be noted that Proposition \ref{main3} does not establish that (VS) admits a unique solution. Instead, uniqueness is a by-product of the verification procedure in Section \ref{section_verification}. The proof of Proposition \ref{main3} consists of three steps, which we develop in Sections \ref{subsection_GeneralODE}--\ref{cons}.

\subsection{An ODE for the Free Boundary} \label{subsection_GeneralODE}

The first step of the proof consists in showing that, if a solution $(\underline m, b, W)$ to (VS) exists, then the function $b$ describing the free boundary satisfies an ODE.
Specifically, let $E: \mathcal J_E \to \mathbb R$ be the vector field defined over $\mathcal J_E \equiv \{(x,m) \in \mathcal J : x \geq x_R\} \setminus \mathcal D$ by
\begin{align}
E(x,m) \equiv \;& \frac{H(m)\sigma^2(x)}{ 2L(x)h_2(x)} \notag \allowdisplaybreaks
\\
& \!\! \times\! \left\{  \frac{\gamma S'(x)}{D(x,m)}\,[R(x)h_2(m)  - G(m ) h_2(x)] + R'(x)h_2(x)- R(x) h'_2(x)\right\} \!, \label{E}
\end{align}
where
\begin{align}
L(x)  &\equiv {\cal L} R(x) -r R(x), \label{L}
\\
D(x,m)&\equiv h_1(m)h_2(x) -h_1(x)h_2(m), \label{D} \allowdisplaybreaks
\\
H (m) &\equiv \frac{f(m)}{1 - F(m)}. \label{H}
\end{align}
Notice that $H$ corresponds to the \textit{breakthrough rate}, that is, $H(m) \, \mathrm dm$ is the probability that a breakthrough occurs over $(m , m+ \mathrm dm]$ conditional on no breakthrough occurring over $(\alpha,m]$. Observe that $E$ is well-defined because $L < 0$ over $[x_R, \beta)$ by A3 along with the fact that $x _R > x_0$, and because $D > 0$ over $\mathcal J_E$ as $h_1$ is strictly increasing and $h_2$ is strictly decreasing. Moreover, $E$ is continuous and, as $h_1$ and $h_2$ are $\mathcal C^2$, $G$ is $\mathcal C^1$, and $H$ is locally Lipschitz by A8, $E$ is locally Lipschitz in its second argument. However, because $D$ vanishes over $\mathcal D$, $E$ cannot be continuously extended to the closure $\overline{\mathcal J_E } \equiv \{(x,m) \in \mathcal J : x \geq x_R\}$ of $\mathcal J_E$, as the following lemma shows.

\begin{lem} \label{explodediadlem}
For each $m \in [x_R, \beta),$
\begin{eqnarray} \label{explodediad}
\lim_{(x,m'), x<m' \to (m,m)} E(x,m')=\infty.
\end{eqnarray}
\end{lem}

\vskip 1.5mm

We are now ready to state the central result of this section.

\begin{lem} \label{odeforb}
If $(\underline m,b,W)$ is a solution to (VS) such that $b \in \mathcal C^1([\underline m, \beta))$ and $W \in \mathcal V,$ then $b$ satisfies the ODE
\begin{align}
b'(m) &= E(b (m),m), \quad m \in [\underline m, \beta), \label{edocase1a}
\\
b(\underline m)  & =  x_R.   \label{edocase1ain}
\end{align}
\end{lem}

\vskip 1.5mm

The proofs of Lemmas \ref{explodediadlem}--\ref{odeforb} are provided in Appendix A.

In light of \eqref{lalala}--\eqref{tc}, Lemma \ref{odeforb} narrows down the set of candidates for the free boundary to the set of solutions to (\ref{edocase1a})--(\ref{edocase1ain}) that satisfy $x_R \leq b(m) < m$ for all $m \in [\underline m,\beta)$ and $\lim_{m \to \beta^-} b(m) = \beta$ for an appropriate choice of the endpoint $\underline m$.

\subsection{The Set of Candidates for the Free Boundary}

The second step of the proof consists in showing that the set of candidates for the free boundary is nonempty. As the vector field $E$ in (\ref{edocase1a}) cannot be continuously extended to $\overline{\mathcal J_E}$ by Lemma \ref{explodediadlem}, the existence of such a candidate does not immediately follow from standard results and requires a specific analysis. Our approach relies on the following lemma, whose geometrical interpretation is that the vector field $E$ points to the left above a strictly increasing $\mathcal C^1$ curve.

\begin{lem} \label{mx}
There exists a strictly increasing $\mathcal C^1$ mapping $x \mapsto m_x$ over $[x_R, \beta)$ such that
\begin{itemize}

\item[(i)]

for each $x \in [x_R,\beta),$ $m_x > x;$

\item[(ii)]

for each $(x,m) \in \mathcal J_E,$ $E(x,m_x)=0$ and $E(x,m) \lessgtr 0 $ if $m \gtrless m_x;$

\item[(iii)]

$\lim_{x \to \beta^-} m_x =\beta$.

\end{itemize}
\end{lem}

\noindent \textbf{Proof.} We use a change of variable introduced by Dayanik and Karatzas (2003). For each $x \in \mathcal I$, define $\zeta(x) \equiv \frac{h_1(x)}{h_2(x)}$, which is strictly increasing in $x$ and maps $\mathcal I$ onto $(0,\infty)$. For any function $g : \mathcal I \to \RR$, define the function $\hat{g}$ by
\begin{eqnarray} \label{changeofvar}
\hat g (y) \equiv \frac{g}{h_2} \circ \zeta^{-1} (y), \quad y \in(0,\infty).
\end{eqnarray}
A direct computation (De Angelis, Ferrari, and Moriarty (2018, Appendix A.1)) shows that, if $g$ is twice differentiable at $x\in \mathcal I$, then $(\CL g-rg)(x)$ has the same sign as $\hat{g}''(\zeta(x))$. Hence A3 and $x_R > x_0$ imply that $\hat{R}''<0$ over $[\zeta(x_R), \infty)$, and A7 and $G$ is $\mathcal C^1$ imply that $\hat G$ is concave over $(0,\infty)$. Moreover, $\hat G >0$ over $(0, \infty)$ as $G > 0$ over $\mathcal I$, and (\ref{gp'}) implies
\begin{eqnarray} \label{gp3}
\lim_{y \to \infty} \frac{\hat{G}(y)}{y}=0.
\end{eqnarray}
Finally, because $R  > 0$ over $[x_R, \beta)$ by Lemma \ref{lem_sign}, we have $\hat R > 0$ over $[\zeta(x_R), \infty)$; because $\hat{R}$ is strictly concave over $[\zeta(x_R), \infty)$, it must then be that $\hat{R}' >0$ over $[\zeta(x_R), \infty)$.

We now use (\ref{changeofvar}) to obtain a more compact expression for $E(x,m)$. By (\ref{wronskian}) and the definition of $\zeta$, we have
\begin{eqnarray*}
\zeta'(x) = \frac{h_1'(x) h_2(x) -h_1(x) h_2'(x)}{[h_2(x)]^2} = {\gamma S'(x) \over [h_2(x)]^2},
\end{eqnarray*}
for all $x \in \mathcal I$, and thus $(\zeta^{-1})'(\zeta(x)) = {[h_2(x)]^2 \over \gamma S'(x)}$. By (\ref{changeofvar}), this implies
\begin{eqnarray}
\hat R'(\zeta (x)) = \! \left(R\over h_2 \right)'\!(x)(\zeta^{-1})'(\zeta(x))={R'(x) h_2(x) - R(x) h_2'(x) \over \gamma S'(x)}. \label{23:42}
\end{eqnarray}
Thus, for each $(x,m) \in \mathcal J_E$, we have
\begin{align}
E(x,m)&=  \frac{H(m)\sigma^2(x)\gamma S'(x)}{ 2L(x)h_2(x)}\!\left[ \frac{R(x)h_2(m)  - G(m ) h_2(x)}{h_1(m)h_2(x) -h_1(x)h_2(m)} + {R'(x)h_2(x)- R(x) h'_2(x) \over \gamma S'(x)}\right]  \notag
\\
&=  \frac{H(m)\sigma^2(x)\gamma S'(x)}{ 2L(x)h_2(x)}\!\left\{ \frac{1}{\zeta(m) - \zeta(x)} \! \left[{R(x)\over h_2(x)}  - {G(m ) \over h_2(m)} \right] + \hat R'(\zeta (x))\right\} \notag \allowdisplaybreaks
\\
&=  -\frac{H(m)\sigma^2(x) \gamma S'(x)}{ 2L(x)h_2(x)} \!\left[ \frac{\hat G(\zeta(m)) - \hat{R} (\zeta(x))}{\zeta(m)-\zeta(x)} - \hat{R}'(\zeta(x)) \right]\hskip -1mm,  \label{Ereecrit}
\end{align}
where the first equality follows from (\ref{E}) and (\ref{D}), the second inequality follows from (\ref{23:42}) and the definition of $\zeta$, and the third equality follows from (\ref{changeofvar}).

Because $L(x) <0$ for all $x \geq x_R$, the upshot from (\ref{Ereecrit}) is that, for each $(x,m) \in \mathcal J_E$,
\begin{eqnarray}
\mathrm{sgn}\, E(x,m) = \mathrm{sgn}\hskip -0.3mm  \left[ \frac{\hat G(\zeta(m)) - \hat{R} (\zeta(x))}{\zeta(m)-\zeta(x)} - \hat{R}'(\zeta(x))\right] \hskip -1mm. \label{Sking}
\end{eqnarray}
Using the notation $z=\zeta(m)$ and $y=\zeta(x)$, consider for all $z>y>0$ the quantity
\begin{eqnarray}
\eta(z,y) \equiv \frac{\hat G(z) - \hat{R} (y)}{z-y} - \hat{R}'(y). \label{eta}
\end{eqnarray}
For each $y\geq \zeta(x_R)$, we have $\lim_{z \to y^+} \eta(z,y)= \infty$ and $\lim_{z \to \infty} \eta(z,y) <0$ by (\ref{gp3}) along with the fact that $\hat R'(y) >0$. Moreover, for each $z >y$,
\begin{eqnarray}
{\partial \eta \over \partial z}\, (z,y) = \frac{\hat G'(z)(z-y) -[\hat G(z)- \hat R (y)] }{(z-y)^2} < \frac{\hat G'(z)(z-y) -[\hat G(z)- \hat G(y)]}{(z-y)^2} \leq 0, \label{etaz}
\end{eqnarray}
where the first inequality follows from $\hat G(y) > \hat R(y)$, and the second inequality follows from the concavity of $\hat G$. This shows that, for each $y \geq \zeta(x_R)$, there exists a unique $z_y > y$ such that $\eta(z_y,y)=0$, $\eta(z,y)<0$ for all $z>z_y$, and $\eta(z,y)>0$ for all $z \in (y, z_y)$. By construction, $\lim_{y \to \infty} z_y=\infty$. To conclude the proof, observe from (\ref{eta}) and $\eta(z_y,y)=0$ that
\begin{eqnarray}
{\partial \eta \over \partial y}\, (z_y,y) = \frac{-\hat R'(y)(z_y-y) + \hat G(z_y) - \hat{R} (y) }{(z_y-y)^2} - \hat R''(y) = - \hat R''(y)>0 \label{etay}
\end{eqnarray}
for all $y \in [\zeta(x_R),\infty)$. Applying the implicit function theorem, we obtain from (\ref{etaz})--(\ref{etay}) that the mapping $y \mapsto z_y$ is $\mathcal C^1$, with ${\mathrm dz_y \over \mathrm d y} >0$ over $[\zeta(x_R), \infty)$. Letting $m_x \equiv \zeta^{-1}(z_{ \zeta(x) })$ for all $x \in [x_R, \beta)$ and recalling that $ \zeta $ is $\mathcal C^2$ and strictly increasing, it is straightforward to verify that the mapping $x \mapsto m_x$ is $\mathcal C^1$ and strictly increasing over $[x_R, \beta)$, and that it satisfies (i)--(iii). The result follows. \hfill $\blacksquare$

\bigskip

It follows from Lemmas \ref{odeforb}--\ref{mx} along with the limit condition $\lim_{m \to \beta^-}b(m)= \beta$ that, in the space $(x,m)$, the free boundary $x=b(m)$ must lie strictly below the locus of points $m= m_x$, at which $E$ vanishes, and strictly above $\mathcal D$, at which $E$ explodes; this situation is illustrated in Figure 2. In particular, the endpoint $\underline m$ must belong to the interval $(x_R , m_{x_R})$. Notice that it is crucial for this argument that the mapping $x \mapsto m_x$ be $\mathcal C^1$. We denote by $\mathcal J_E^+ \equiv \{(x,m) \in \mathcal J_E : m_x > m > x\}$ the corresponding domain.
\noindent
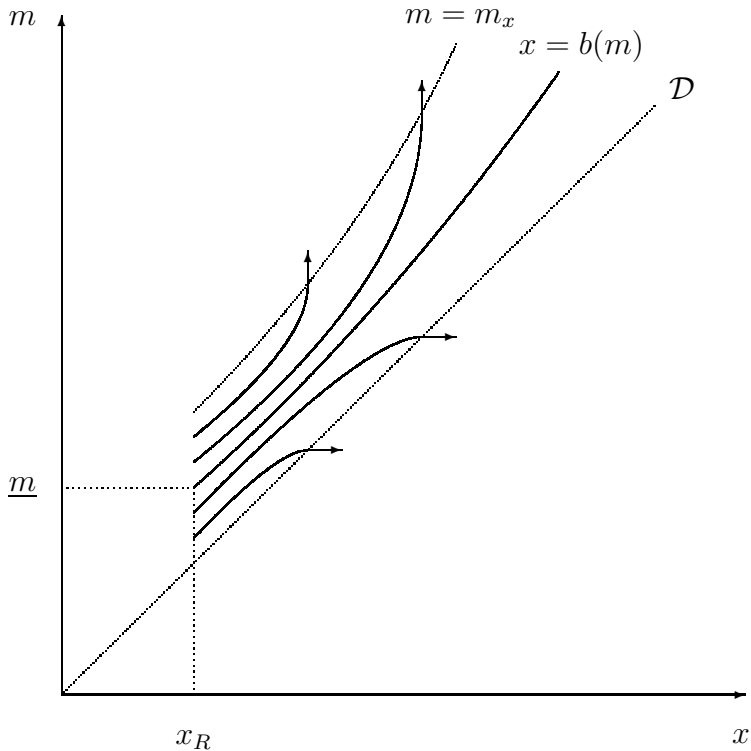
\begin{figure}
{\normalsize \noindent \unitlength=1.50mm \special{em:linewidth
0.4pt} \linethickness{0.4pt} \ \ \ \ \ \ \ \ \
\begin{picture}(140.00,70.00)
\put(18.4,10){\vector(1,0){60}}
\put(18.4,10){\vector(0,1){60}}
\put(78,6.28){\makebox(0,0)[cc]{$x$}}
\put(15,69.8){\makebox(0,0)[cc]{$m$}}
\bezier{242}(18.4,10)(44.4,36)(70.4,62)
\put(72.8,63.6){\makebox(0,0)[cc]{$\mathcal D$}}
\bezier{150}(30,35)(45,50.1)(53,67.5)
\put(53.4,69.8){\makebox(0,0)[cc]{$m=m_x$}}


\put(40,46.2){\vector(0,1){3}}
\bezier{1000}(30,32.77)(40,41)(40,46.2)

\put(50,61.3){\vector(0,1){3}}
\bezier{1000}(30,30.54)(50,47)(50,61.3)


\put(40,31.6){\vector(1,0){3}}
\bezier{1000}(30,23.83)(37.3,31.6)(40,31.6)

\put(50,41.6){\vector(1,0){3}}
\bezier{1000}(30,26.06)(45,41.6)(50,41.6)


\bezier{1000}(30,28.29)(47,43.5)(62,65)
\put(64,67.3){\makebox(0,0)[cc]{$x= b(m)$}}
\put(15,28.29){\makebox(0,0)[cc]{$\underline m$}}
\put(30,6){\makebox(0,0)[cc]{$x_R$}}

\bezier{36}(30,28.29)(30,19.145)(30,10)
\bezier{23}(30,28.29)(24.34,28.29)(18.4,28.29)

\end{picture}}\vskip -7mm
\caption{The vector field $E$ and the free boundary $x = b(m)$.}
\end{figure}

The set of candidates for the free boundary $b$ can then be formally described as follows. For each $m_0 \in (x_R, m_{x_R})$, let $b_{m_0}$ be the maximal solution to the ODE $b'(m) = E(b(m),m)$ with initial condition $b(m_0) = x_R$ that satisfies $(b_{m_0 }(m),m) \in \mathcal J_E^+$ for all $m$ in a nonempty maximal interval $(m_0, \overline m_{m_0})$. Because the vector field $E$ is continuous and locally Lipschitz in its second argument and $b'_{m_0}(m_0) > 0$ as $m_0 < m_{x_R}$, the existence and uniqueness of $b_{m_0}$ is guaranteed for all $m_0 \in (x_R, m_{x_R})$ by the Cauchy--Lipschitz theorem. For any such $m_0$, we will say that $b_{m_0}$ is a \textit{candidate for the free boundary} if and only if $\overline m_{m_0} = \beta$. Observe that we then have $m_{b_{m_0}(m)} > m$ for all $m\in [m_0, \beta)$; hence, by Lemma \ref{mx}(iii), the limit condition  $\lim_{m \to \beta^-} b_{m_0} (m)=\beta$ is automatically satisfied, which justifies the terminology. Moreover, $b_{m_0}$ is strictly increasing as $\mathcal J_E^+= E^{-1}( (0, \infty))$ by Lemma \ref{mx}(ii). As the following lemma shows, the fact that the vector field $E$ is outward-pointing at the boundary of $\mathcal J_E^+$ guarantees the existence of such candidates.

\begin{lem} \label{b}
There exists a nonempty compact interval $\mathcal I^0 \subset (x_R, m_{x_R})$ such that $b_{m_0}$ is a candidate for the free boundary if and only if $m_0 \in \mathcal I^0$.
\end{lem}

\noindent \textbf{Proof.} Fix some $m_0 \in (x_R,m_{x_R})$ and the corresponding maximal solution $b_{m_0}$. If $\overline m_{m_0} <\beta$, then we either have $b_{m_0}(\overline m_{m_0}^-)=\overline m_{m_0}$ and $b_{m_0}$ cuts $\partial ^- \mathcal J_E ^+\equiv \{(x,m) \in \mathcal D: x \geq x_R\}$ at $\overline m_{m_0}$, or $b_{m_0}(\overline m_{m_0}^-) = m_{b_{m_0}(\overline m_{m_0}^-)}$ and $b_{m_0}$ cuts $\partial ^ + \mathcal J_E^+ \equiv \{(x,m) \in \mathcal J_E: m=m_x\}$ at $\overline m_{m_0}$. Define also the degenerate solutions $b_{x_R} \equiv \{(x_R,x_R)\}$ and $b_{m_{x_R}}= \{(x_R,m_{x_R})\}$, which correspond to the limit cases of these two situations, respectively. Therefore, we can partition the interval $[x_R,m_{x_R}]$ into three pieces:
\begin{enumerate}

\item

$\mathcal I^-$, the set of $m_0$ such that $b_{m_0}$ cuts $\partial ^-\mathcal J_E^+$;

\item

$\mathcal I^0$, the set of $m_0$ such that $\overline m_{m_0}=\beta$;

\item

$\mathcal I^+$, the set of $m_0$ such that $b_{m_0}$ cuts $\partial ^+ \mathcal J_E^+$.

\end{enumerate}
Clearly $\mathcal I^- \neq \emptyset $ as $x_R \in \mathcal I^-$, $\mathcal I^+ \neq \emptyset $ as $m_{x_R} \in \mathcal I^+$, and $\mathcal I^- \cap \mathcal I^+ = \emptyset$ as $m_x >x$ for all $x \in [x_R, \beta)$. From the non-crossing property of the solutions to the ODE $b'(m) = E(b(m),m)$, $\mathcal I^-$ and $\mathcal I^+$ are intervals. If $\mathcal I^-$ and $\mathcal I^+$ are relatively open in $[x_R, m_{x_R}]$, then, because $[x_R, m_{x_R}]$ cannot be the union of two disjoints open intervals, $\mathcal I^0$ must be a nonempty closed interval in $(x_R, m_{x_R})$, which concludes the proof. The fact that $\mathcal I^-$ and $\mathcal I^+$ are relatively open in $[x_R, m_{x_R}]$ follows directly along the lines of the proof of Theorem 1 in Bobtcheff, Bolte, and Mariotti (2017), to which we refer for details. Specifically, that $\mathcal I^-$ is relatively open in $[x_R, m_{x_R}]$ is a consequence of the fact that, by Lemma \ref{explodediadlem}, the vector field $E$ explodes over $\partial ^-\mathcal J_E^+$, whose slope is 1; and that $\mathcal I^+$ is relatively open in $[x_R, m_{x_R}]$ is a consequence of the fact that, by Lemma \ref{mx}, the vector field $E$ vanishes over $\partial ^+ \mathcal J_E^+$, whose slope is locally bounded as the mapping $x \mapsto m_x$ is $\mathcal C^1$. The result follows. \hfill $\blacksquare$

\bigskip

The proof essentially follows the retraction principle of Wa$\dot{\mathrm z}$ewski (1947), see for instance Hartman (1964, Chapter X, Theorem 2.1), with slight adjustments owing to the fact that $E$ cannot be continuously extended to $\overline {\mathcal J_E}$. Observe that Lemma \ref{b} shows the existence but not the uniqueness of a candidate for the free boundary: for all we know at this stage, the set $\mathcal I^0$ may not be reduced to a singleton.

\subsection{A Solution to the Variational System} \label{cons}

The third and final step of the proof consists in showing that to every candidate $b: [\underline m, \beta) \to [x_R, \beta)$ for the free boundary corresponds a solution $(\underline m, b, W)$ to (VS). Specifically, because any such $b \in \mathcal C^1([\underline m, \beta))$ is strictly increasing and satisfies (\ref{ic})--(\ref{tc}), all we need to do in order to complete the proof of Proposition \ref{main3} is to exhibit a function $W \in \mathcal V$ that satisfies (\ref{vs1})--(\ref{vs5}) for $\mathcal S \equiv \{ (x,m) \in \mathcal J : \underline m \leq m <\beta  \text{ and } x_R \leq x \leq  b(m)\}$ and $\mathcal C \equiv \mathcal J \setminus \mathcal S$. Our construction of $W$ is explicit given $b$ and consists of four steps.

\subparagraph{Step 1}

Consider first the \pagebreak points $(x,m)\in \mathcal S$. In this region, we let $W$ be equal to the payoff from stopping $(X,M)$ immediately,
\begin{eqnarray}\label{Z1}
W (x,m) \equiv [1-F(m)]R(x), \quad (x,m)\in \mathcal S,
\end{eqnarray}
in line with (\ref{vs2}).

\subparagraph{Step 2}

Consider next the points $(x,m) \in \mathcal C$ to the left of $\mathcal S$, such that $(x,m) \in (\alpha, x_R) \times [\underline  m, \beta )$. In this region, we let $W$ be equal, up to multiplication by $1- F(m)$, to the stand- alone value function $V_R$,
\begin{eqnarray} \label{Z2}
{W}(x,m) \equiv [1 - F(m)] \, \frac{h_1(x)}{h_1(x_R)} \, R(x_R), \quad (x,m) \in  (\alpha, x_R) \times [\underline  m, \beta ),
\end{eqnarray}
which satisfies (\ref{vs1}) by definition of $h_1$ and (\ref{vs5}) by (\ref{hboundaries}), and pastes continuously with (\ref{Z1}).

\subparagraph{Step 3}

Consider now the points $(x,m)\in \mathcal C$ to the right of $\mathcal S$, such that $m \in [\underline m, \beta)$ and $x \in  (b(m), m]$. In this region, we let $W$ be equal to the solution to \eqref{vs1} that pastes continuously with (\ref{Z1}) at the free boundary $x=b(m)$, and that satisfies the smooth-fit condition (\ref{vs3}). As shown in the proof of Lemma \ref{odeforb}, this leads to
\begin{eqnarray}\label{Z3}
W(x,m) \equiv A(m)h_1(x)+B(m)h_2(x),\quad m \in [\underline m, \beta) \text{ and } x \in (b(m),m],
\end{eqnarray}
where
\begin{align}
A(m)  &= \frac{1-F(m)}{\gamma S'(b(m))} \,[R'(b(m))h_2( b(m))- R(b(m)) h'_2(b(m))], \label{exprA}
\\
B(m) &= - \frac{1-F(m)}{\gamma S'(b(m))} \,[R'(b(m))h_1( b(m))- R(b(m)) h'_1(b(m))], \label{exprB}
\end{align}
and the ODE (\ref{edocase1a}) for $b$ guarantees that the Neumann condition (\ref{vs4}) is satisfied. Notice that $B(\underline m)=0$ as $b(\underline m)=x_R$ by construction and $R'(x_R)=\frac{h'_1(x_R)}{h_1(x_R)}\,R(x_R)$ by the smooth-fit property for $V_R$.

\subparagraph{Step 4}

Consider finally the remaining points $(x,m)\in \mathcal C$, such that $m \in (\alpha, \underline m)$ and $x \in  (\alpha, m]$. In this region, we let $W$ be equal to the solution to \eqref{vs1} that satisfies the boundary condition (\ref{vs5}), the Neumann condition (\ref{vs4}), and that pastes continuously with the solution constructed so far. By (\ref{hboundaries}), the boundary condition leads to
\begin{eqnarray}\label{Z4}
W(x,m)\equiv C(m)h_1(x),\quad m \in (\alpha, \underline m) \mbox{ and } x \in  (\alpha, m],
\end{eqnarray}
for some function $C \in \mathcal C^1((\alpha,\underline m))$. Because $B(\underline m) = 0$, continuous pasting requires
\begin{eqnarray*}
C(\underline m)= A(\underline m) = {1-F(\underline m) \over h_1(x_R)} \,R(x_R)
\end{eqnarray*}
by (\ref{exprA}), using (\ref{wronskian}) along with the smooth-fit property for $V_R$. The Neumann condition is satisfied if and only if
\begin{eqnarray*}
C'(m)= -\frac{f(m)}{h_1(m)}\, G(m),
\end{eqnarray*}
from which we conclude that
\begin{eqnarray}\label{exprAbelow}
C(m) \equiv {1-F(\underline m) \over h_1(x_R)} \,R(x_R) +\int_m^{\underline m}\frac{f(y)}{h_1(y)}\,G(y)\, \mathrm dy.
\end{eqnarray}
This completes the construction of a function $W$ satisfying (\ref{vs1})--(\ref{vs5}) and, hence, of a solution $(\underline m, b, W)$ to (VS). It should be noted that, given $b$, $W$ is the unique solution to (\ref{vs1})--(\ref{vs5}) that is continuous over $\mathcal J$. Together with Lemma \ref{b}, which characterizes the set of candidates for the free boundary $b$, this yields a complete characterization of the solutions to (VS).

\vskip 3mm

There only remains to check that the function $W$ defined by (\ref{Z1})--(\ref{exprAbelow}) satisfies the regularity conditions required in Proposition \ref{main3}.

\begin{lem} \label{WinF}
$W\in \mathcal V$.
\end{lem}

The proof of Lemma \ref{WinF} is provided in Appendix A.

To conclude this section, we note that $W$ is strictly positive over $\mathcal J$, which reflects that the functions $A$, $B$, and $C$ are strictly positive over $[\underline m,\beta)$, $(\underline m,\beta)$, and $(\alpha, \underline m]$, respectively.

\begin{cor} \label{corosigne}
For each $(x,m)\in \mathcal J,$ $W(x,m) >0$.
\end{cor}

\section{Verification} \label{section_verification}

Let us now fix a candidate $b$ for the free boundary, which by Lemmas \ref{odeforb} and \ref{b} satisfies the ODE \eqref{edocase1a}--\eqref{edocase1ain} for some endpoint $\underline m \in \mathcal I^0$ and is strictly increasing, and let $(\underline m,b, W)$ be the corresponding solution to (VS) constructed in Section \ref{cons}. Denote by
\begin{eqnarray*}
\tau_b \equiv \inf\hskip 0.5mm \{ t\geq 0 : (X_t,M_t)\in \mathcal S\}
\end{eqnarray*}
the hitting time of the stopping region $\mathcal S = \{ (x,m) \in \mathcal J : \underline m \leq m <\beta  \text{ and } x_R \leq x \leq  b(m)\}$ associated to $b$, which need not be $\mathbf P_{x,m}$-almost surely finite, and denote by
\begin{eqnarray}\label{defW}
V_b(x,m) \equiv \mathbf E_{x,m} \! \left[[1-F(M_{\tau_b})]\, \mathrm e^{-r \tau_b } R(X_{\tau_b }) +\int_{0}^{\tau_b }\mathrm e^{-r t}G(M_t)f(M_t)\, \mathrm dM_t  \right]
\end{eqnarray}
the value for the DM of stopping $(X,M)$ at $\tau_b$, starting from $(x,m)$.

So far, our analysis has thus led us to construct three functions over $\mathcal J$: the value function $V$ of problem (\ref{mainproblem}), the function $W$ analytically characterized by (\ref{Z1})--(\ref{exprAbelow}) as part of a solution $(\underline m, b, W)$ to (VS) given the candidate $b$ for the free boundary, and the function $V_b$ defined by (\ref{defW}) given the stopping time $\tau_b$. The central result of this section is that these three functions coincide.

\begin{pro} \label{=}
$V = W = V_b$.
\end{pro}

Proposition \ref{=} concludes the proof of Theorem \ref{maintheo}. As a by-product, it entails that the interval $\mathcal I^0$ of possible endpoints for a candidate for the free boundary is reduced to a point, which in turn implies that there is a single such candidate. Indeed, if there were two different candidates $b_1$ and $b_2$ for the free boundary corresponding to different endpoints $\underline m\,\!_1, \underline m\!\,_2 \in \mathcal I^0$, then we would have two different solutions $(\underline m\,\!_1, b_1, W_1)$ and $(\underline m\,\!_1, b_1, W_2)$ to (VS), contradicting Proposition \ref{=} given that the value function $V$ of problem (\ref{mainproblem}) is uniquely defined. Therefore, (VS) has a unique solution.

The proof of Proposition \ref{=} consists of three steps, which we develop in Sections \ref{Garance}--\ref{Sreemati}.

\subsection{Two Useful Bounds} \label{Garance}

The first step of the proof consists in providing appropriate bounds for the function $W$ defined by (\ref{Z1})--(\ref{exprAbelow}). We first show that $W$ is bounded below by the DM's payoff from stopping $(X,M)$ immediately.

\begin{lem} \label{lem_above}
For each $(x,m) \in \mathcal J,$
\begin{eqnarray}
W(x,m) \geq [1-F(m)]R(x), \label{geqobstacle}
\end{eqnarray}
and this inequality is strict if $(x,m) \in \mathcal C$.
\end{lem}

\noindent \textbf{Proof.} By \eqref{Z1}, (\ref{geqobstacle}) holds as an equality for $(x,m) \in \mathcal S$. We thus only need to prove that (\ref{geqobstacle}) holds as a strict inequality for $(x,m) \in \mathcal C$. We consider two cases in turn.

\subparagraph{Case 1}

Suppose first that $(x,m) \in \mathcal C$ is such that $m \in (\alpha, \underline m)$, so that $W(x,m)$ is given by \eqref{Z4}--\eqref{exprAbelow} and (\ref{geqobstacle}) is equivalent to
\begin{eqnarray*}
\frac{C(m)}{1-F(m)} \geq \frac{R(x)}{h_1(x)}.
\end{eqnarray*}
Now, observe that the optimality of the stopping threshold $x_R$ for problem (\ref{sa}) implies that, for each $x \in \mathcal I$,
\begin{eqnarray*}
\frac{R(x_R)}{h_1(x_R)} \geq \frac{R(x)}{h_1(x)}.
\end{eqnarray*}
Thus it is sufficient to show that, for each $m \in (\alpha, \underline m)$,
\begin{eqnarray*}
\frac{C(m)}{1-F(m)} > \frac{R(x_R)}{h_1(x_R)}
\end{eqnarray*}
or, equivalently, by (\ref{exprAbelow}),
\begin{eqnarray}\label{eq_sufficient_cond}
\int_m^{\underline m} \frac{f(y)}{h_1(x_R)}\! \left[ \frac{h_1(x_R)}{h_1(y)}\, G(y)-R(x_R)\right] \! \mathrm dy > 0.
\end{eqnarray}
We now prove that the integrand in \eqref{eq_sufficient_cond} is strictly positive for all $y \in (\alpha,\underline m)$, which concludes the discussion of this case. For $y \in(\alpha, x_R)$, we have
\begin{eqnarray*}
\frac{h_1(y)}{h_1(x_R)}\, G(x_R)=\mathbf E_y\hskip 0.3mm[\mathrm e^{-r \tau(x_R)}G(X_{\tau(x_R)})]= G(y)+ \mathbf E_y \! \left[\int_0^{\tau(x_R)}\mathrm e^{-rt}(\CL G-rG)(X_t)\, \mathrm dt \right] \! \leq G(y)
\end{eqnarray*}
by It\^o's lemma, where the inequality follows from A7. Hence, by A6,
\begin{eqnarray*}
\frac{h_1(x_R)}{h_1(y)}\, G(y) - R(x_R) \geq G(x_R) - R(x_R)  > 0,
\end{eqnarray*}
as desired. For $y \in[x_R,\underline m)$, we similarly have
\begin{eqnarray*}
\frac{h_1(x_R)}{h_1(y)}\, G(y)= \mathbf E_{x_R}\hskip 0.3mm[\mathrm e^{-r \tau(y)}G(X_{\tau(y)})]= G(x_R)+\mathbf E_{x_R}\!\left[\int_0^{\tau(y)} \mathrm e^{-rt}(\CL G- rG)(X_ t) \, \mathrm dt\right]\hskip -1mm.
\end{eqnarray*}
Because the mapping $y \mapsto \tau(y)$ is $\mathbf P_{x_R}$-almost surely strictly increasing, it follows from this and A7 that the mapping $y \mapsto \frac{h_1(x_R)}{h_1(y)}\, G(y)$ is nonincreasing over $[x_R, \underline m)$. To conclude, we thus only need to check that
\begin{eqnarray}
\frac{h_1(x_R)}{h_1(\underline m)} \, G(\underline m)> R(x_R). \label{comparpayoffs}
\end{eqnarray}
To see this, recall that, because $\underline m \in \mathcal I^0 \subset (x_R,m_{x_R})$ by Lemma \ref{b}, we have $E(x_R, \underline m) >0$ by Lemma \ref{mx}. Hence
\begin{align*}
1  &= \mathrm {sgn} \left\{  \frac{\gamma S'(x_R)}{D(x_R,\underline m)}\,[G(\underline m) h_2(x_R)- R(x_R)h_2( \underline m) ] +R(x_R) h'_2(x_R) - R'(x_R)h_2(x_R) \right\}
\\
&= \mathrm {sgn} \left\{  \frac{\gamma S'(x_R)}{D(x_R,\underline m)}\,[G(\underline m) h_2(x_R)- R(x_R)h_2( \underline m) ] +R(x_R) \! \left[h'_2(x_R) - {h'_1(x_R) \over h_1(x_R)}\,h_2(x_R)\right] \right\}
\\
&= \mathrm {sgn} \left\{\frac{1}{D(x_R,\underline m)}\,[G(\underline m) h_2(x_R)- R(x_R)h_2( \underline m) ] - {1 \over h_1(x_R)}\, R(x_R) \right\}
\\
&= \mathrm {sgn} \left[\frac{h_1(x_R)}{h_1(\underline m)} \, G(\underline m)- R(x_R) \right]\hskip -1mm,
\end{align*}
where the first equality follows from (\ref{E}) along with the fact that $L(x_R) <0$, the second equality follows from the smooth-fit condition $R'(x_R)=\frac{h'_1(x_R) }{h_1(x_R)} \, R(x_R)$ for problem (\ref{sa}), the third inequality follows from (\ref{wronskian}), and the fourth equality follows from (\ref{D}) along with the fact that $D(x_R, \underline m) >0$. This proves (\ref{comparpayoffs}), which concludes the discussion of this case.

\subparagraph{Case 2}

Suppose next that $(x,m) \in \mathcal C$ is such that $m \geq \underline m$. For $x \in(\alpha, x_R)$, we have $W(x,m) = [1 - F(m)] V_R(x)$ by (\ref{Z2}), and the result follows from noticing that $V_R > R$ over $(\alpha,x_R)$. For $x \in (b(m), m]$ and, hence, $x >x_R$, we once more rely on the change of variable (\ref{changeofvar}). As pointed out in the proof of Lemma \ref{mx}, A3 and $x_R > x_0$ imply that $\hat{R}''<0$ over $[\zeta(x_R), \infty)$; thus $\hat R$ is strictly concave over $[\zeta(x_R), \infty)$. Now, let us fix some $m \geq \underline m$ and, for each $x \geq b(m)$, let $u_m(x) \equiv W(x,m)$. By (\ref{changeofvar}) and the continuous-fit condition $u_m(b(m))=[1-F(m)]R(b(m))$, we have
\begin{eqnarray}
\hat u_m(\zeta(b(m)))=[1-F(m)]\hat R(\zeta(b(m))). \label{And}
\end{eqnarray}
Next, by (\ref{changeofvar}) and (\ref{Z3}), we have
\begin{eqnarray}
\hat u _m (y)= A(m)y+ B(m) \label{Elo}
\end{eqnarray}
for all $y \geq \zeta(b(m))$. Finally, by (\ref{23:42}) and (\ref{exprA}), we have
\begin{eqnarray}
[1-F(m)]\hat R'(\zeta(b(m))) = [1-F(m)]\,{R'(b(m)) h_2(b(m)) - R(b(m)) h_2'(b(m)) \over \gamma S'(b(m))} = A(m). \label{Ale}
\end{eqnarray}
Taken together, (\ref{And})--(\ref{Ale}) imply that the affine mapping $y \mapsto \hat u_m(y)$ is tangent to the mapping $y \mapsto [1-F(m)]\hat{R}(y)$ at $\zeta(b(m))$. As the latter is strictly concave over $[\zeta(x_R), \infty)$ and, hence, over $[\zeta(b(m)), \infty)$, we obtain that
\begin{eqnarray*}
\hat u_m(y) > [1-F(m)]\hat R(y)
\end{eqnarray*}
for all $y \in (\zeta(b(m)),\zeta(m)]$, and thus, by (\ref{changeofvar}), that
\begin{eqnarray*}
u_m(x) > [1-F(m)] R(x)
\end{eqnarray*}
for all $x \in (b(m), m]$, which is the desired inequality. The result follows. \hfill $\blacksquare$

\bigskip

It should be noted that the proof of Lemma \ref{lem_above} makes essential use of A7, which reflects that $G$ incorporates itself the solution to an optimal stopping problem. The key inequality (\ref{comparpayoffs}) intuitively expresses that, at $(x_R, \underline m)$, the DM would be ready to wait until $X$ reaches $\underline m$ again if he were certain that a breakthrough would occur at this point, allowing him to obtain a payoff $G(\underline m)$ instead of his payoff $R(x_R)$ from stopping $(X,M)$ immediately; yet, precisely because a breakthrough may occur only later on, it will be optimal for him to invest in the stand-alone technology at $(x_R, \underline m)$.

We next show that $W$ is bounded above by the value the DM could obtain if $G(X)$ were immediately available.

\begin{lem}\label{lem_uniform_integrability}
For each $(x,m) \in \mathcal J,$
\begin{eqnarray}
{W}(x,m) < [1-F(m)] G(x). \label{ineqq}
\end{eqnarray}
\end{lem}

\vskip 1.5mm

\noindent \textbf{Proof.} We consider four cases in turn.

\subparagraph{Case 1}

Suppose first that $(x,m) \in \mathcal S$, so that $W(x,m)$ is given by (\ref{Z1}). Then
\begin{eqnarray*}
W(x,m) = [1 - F(m)] R(x) < [1-F(m)] G(x)
\end{eqnarray*}
by A6.

\subparagraph{Case 2}

Suppose next $(x,m) \in (\alpha,x_R) \times [\underline m, \beta)$, so that $W(x,m)$ is given by (\ref{Z2}). Then
\begin{eqnarray*}
W(x,m)= [1 - F(m)] \, \frac{h_1(x)}{h_1(x_R)} \, R(x_R)  < [1-F(m)] G(x)
\end{eqnarray*}
by (\ref{VR}) and A6.

\subparagraph{Case 3}

Suppose now that $m \in [\underline m, \beta)$ and $x \in  (b(m), m]$, so that $W(x,m)$ is given by (\ref{Z3})--(\ref{exprB}). Then
\begin{align*}
W(x,m) &= A(m)h_1(x)+B(m)h_2(x) \allowdisplaybreaks
\\
&=  [A(m) \zeta(x) + B(m) ] h_2(x)
\\
&= [1-F(m)]\{\hat{R}(\zeta(b(m)))+ \hat{R}'(\zeta(b(m)))[\zeta(x)-\zeta(b(m))]\}h_2(x)
\\
& =  [1-F(m)] \{\hat G (\zeta(x)) - \eta(\zeta(x), \zeta(b(m)))[\zeta(x)-\zeta(b(m))]\}  h_2(x)
\\
& <  [1-F(m)] \hat G (\zeta(x)) h_2(x)
\\
&=  [1-F(m)] G(x),
\end{align*}
where the second equality follows from the definition of $\zeta(x)$, the third equality follows from (\ref{exprA})--(\ref{exprB}) using (\ref{changeofvar}) and (\ref{23:42}), the fourth equality follows from \eqref{eta}, and the inequality follows from noticing that $\zeta (x) > \zeta (b(m))$ and $\eta(\zeta(x), \zeta(b(m))) >0$. To prove this last inequality, let $z_{ \zeta(b(m))} $ be the unique $z$ such that $\eta(z, \zeta (b(m))) =0$ as in the proof of Lemma \ref{mx}; we then have $\eta(z, \zeta(b(m)))>0$ if and only if $z\in (\zeta(b(m)), z_{\zeta(b(m))})$. Because $\zeta(x) > \zeta(b(m))$, we thus only need to show that $\zeta(x) < z_{\zeta(b(m))}$; this follows from the fact that $\zeta(x)\leq \zeta(m)$ along with the observation that (\ref{Sking})--(\ref{eta}) and $E(b(m),m) = b'(m) >0$ imply $\eta(\zeta(m), \zeta(b(m)) >0$ and, hence, $\zeta(m) < z_{\zeta(b(m))}$.

\subparagraph{Case 4}

Suppose finally that $m \in (\alpha, \underline m)$ and $x \in  (\alpha, m]$, so that $W(x,m)$ is given by (\ref{Z4})--(\ref{exprAbelow}). Then
\begin{align*}
W(x,m)   &= \! \left[{1-F(\underline m) \over h_1(x_R)} \,R(x_R) +\int_m^{\underline m}\frac{f(y)}{h_1(y)}\,G(y)\, \mathrm dy\right] \! h_1(x) \allowdisplaybreaks
\\
&< \! \left[{1-F(\underline m) \over h_1(\underline m)} \,G(\underline m) +\int_m^{\underline m}\frac{f(y)}{h_1(y)}\,G(y)\, \mathrm dy\right] \! h_1(x)
\\
&\leq [1- F(m)] G(x),
\end{align*}
where the first inequality follows from \eqref{comparpayoffs}, and the second inequality follows from the fact that the mapping $y\mapsto {G(y) \over h_1(y)}$ is nonincreasing over $[x, \underline m]$, as shown as in the proof of Lemma \ref{lem_above}. The result follows. \hfill $\blacksquare$

\bigskip

It should be noted that the proof of Lemma \ref{lem_uniform_integrability} makes essential use of the fact that the candidate for the free boundary is strictly increasing. Observe that the bounds \eqref{geqobstacle} and \eqref{ineqq} have intuitive interpretations in terms of problem (\ref{mainproblem}): the lower bound (\ref{geqobstacle}) results from letting $\tau \equiv 0$ in (\ref{mainproblem}), while the upper bound \eqref{ineqq} results from (\ref{mainproblem}) upon observing that $G > R$ by A6 and that $(\mathrm e ^{-rt} G(X_t))_{t \geq 0}$ is a supermartingale by A7.

The final result of this section is an immediate consequence of A4, Corollary \ref{corosigne}, and Lemma \ref{lem_uniform_integrability}.

\begin{cor} \label{zoe}
The family $(\mathrm e^{-r\tau} W (X_{\tau},M_{ \tau}))_{\tau \in {\cal T}_{X, M}}$ is uniformly integrable.
\end{cor}

Corollary \ref{zoe} is key to the next step of the verification argument, as we shall now see.

\subsection{The Superharmonicity Property} \label{Gaspard}

The second step of the proof consists in showing that the function $W$ defined by (\ref{Z1})--(\ref{exprAbelow}) is superharmonic.

\begin{lem}
\label{lem_superharmonic}
For all $\tau \in\mathcal T_{x,m}$ and $(x,m) \in \mathcal J,$
\begin{eqnarray}\label{superharmonic}
W(x,m)\geq \mathbf E_{x,m} \! \left[\mathrm e^{-r \tau} W (X_\tau,M_\tau) +\int_{0}^\tau \mathrm e^{-r t}G(M_t)f(M_t) \, \mathrm dM_t \right]\hskip -1mm.
\end{eqnarray}
\end{lem}

\vskip 1.5mm

\noindent \textbf{Proof.} Because $b\in \mathcal C^1([\underline m, \beta))$ is strictly increasing, $(b(M_t))_{t \geq 0}$ is a semimartingale of locally bounded variation. Therefore, we can apply the generalized version of It\^o's lemma due to Peskir (2007, Theorem 4.1).\footnote{Actually, we need a slightly more general version of Peskir's (2007) formula that can easily be obtained by concatenation across the four regions over which $W$ is piecemeal constructed.} We obtain that, for any $\mathbf P_{x,m}$-almost surely finite stopping time $\tau_f \in \mathcal T_{X,M}$,
\begin{align}
W(x,m)  = \;&\mathrm e^{-r \tau_f} W(X_{\tau_f},M_{\tau_f}) + \int_0^{\tau_f} \mathrm e^{-rt} (rW-\mathcal{L} W)(X_t,M_t)1_{\{X_t \notin \{ x_R, b(M_t)\}\} } \, \mathrm dt \nonumber
\\
&\! \!+ \int_0^{\tau_f} \mathrm e^{-rt}f(M_t)G(M_t) \, \mathrm dM_t - \int_0^{\tau_f} \mathrm e^{-rt}\, \frac{\partial W}{\partial x}\,(X_t,M_t) \, \sigma(X_t) \, \mathrm dW_t, \label{itoforW}
\end{align}
$\mathbf P_{x,m}$-almost surely; notice that there is no local-time term in \eqref{itoforW} as $W$ is $\mathcal C^1$ on each horizontal line. Let $([\alpha_n,\beta_n] ) _{n \in \mathbb N}$ be an increasing sequence of compact intervals of $\mathcal I$ such that $\bigcup _{n \in \mathbb N}\, [\alpha_n,\beta_n] = \mathcal I$, and, for each $n \in \mathbb N$, let $\tau_n \equiv \inf \hskip 0.5mm \{ t \geq 0 : X_t \notin [\alpha_n, \beta_n]\}$. Observe that $\tau_n < \infty$ and that $(X_t,M_t)\in [\alpha_n, \beta_n] \times [\alpha_n,\beta_n]$ over $\{t \leq \tau_n\}$, $\mathbf P_{x,m}$-almost surely (Karatzas and Shreve (1991, Chapter 5, Section 5.C)); as $W \in \mathcal V$ and $\sigma$ is continuous, there exists $K_n >0 $ such that $\left|\frac{\partial W}{\partial x} \,(X_t,M_t)\right| \sigma(X_t)\leq K_n$ over $\{t \leq \tau_n\}$, $\mathbf P_{x,m}$-almost surely. It follows that
\begin{eqnarray*}
\mathbf E_{x,m} \! \left[\int_0^{\tau \wedge \tau_n} \mathrm e^{-rt}\,\frac{\partial W}{\partial x} \, (X_t,M_t) \,\sigma(X_t) \, \mathrm dW_t \right] \! = 0.
\end{eqnarray*}
Next, because $\mathcal{L}  W - rW = \mathcal L V_R -rV_R \leq 0$ over $\textrm {int}\, \mathcal S$ and $\mathcal{L}  W - rW = 0$ over $\mathcal C$, we have
\begin{eqnarray*}
\mathbf E_{x,m}\! \left[\int_0^{\tau\wedge \tau_n} \mathrm e^{-rt} (r  W-\mathcal{L}W)(X_t,M_t)1_{\{X_t \notin \{ x_R, b(M_t)\}\} }\, \mathrm dt \right] \! \geq 0.
\end{eqnarray*}
Therefore, applying (\ref{itoforW}) to $\tau_f \equiv \tau \wedge \tau_n$ and taking expectations, we obtain
\begin{eqnarray}\label{eq_ito_dynkin_1}
W(x,m)\geq  \mathbf E_{x,m} \! \left[ \mathrm e^{-r\tau\wedge \tau_n} W(X_{\tau\wedge \tau_n},M_{\tau\wedge \tau_n}) + \int_0^{\tau\wedge \tau_n} \mathrm e^{-rt}G(M_t)f(M_t) \, \mathrm dM_t \right] \hskip -1mm.
\end{eqnarray}
Using Lemma \ref{lemchange-of-time} along with the fact that $G >0$ over $\mathcal I$, we have
\begin{eqnarray*}
\int_0^{\tau\wedge \tau_n} \mathrm e^{-rt}G(M_t)f(M_t) \, \mathrm dM_t \leq \int_0^\infty \mathrm  e^{-rt}G(M_t)f(M_t) \, \mathrm dM_t =\int_m^{\beta} \mathrm e^{-r \tau(y)} G(y)f(y)\, \mathrm dy,
\end{eqnarray*}
and the right-hand side of this inequality has finite expectation as shown in the proof of Lemma \ref{lemchange-of-time}. Because $X$ is a regular diffusion, $\lim _{n \to \infty} \tau_n = \infty$ and, hence, $\lim_{n \to \infty} \tau \wedge \tau _n= \tau$, $\mathbf P_{x,m}$-almost surely. Therefore, we have
\begin{eqnarray}
\lim _{n \to \infty} \mathbf E_{x,m}\! \left[ \int_0^{\tau\wedge \tau_n} \mathrm e^{-rs}G(M_t)f(M_t)\, \mathrm dM_t\right]\! = \mathbf E_{x,m} \! \left[ \int_0^{\tau} \mathrm e^{-rt}G(M_t)f(M_t) \, \mathrm dM_t\right] \hskip -1mm \label{etde1}
\end{eqnarray}
by Lebesgue's dominated convergence theorem. Next, because
\begin{eqnarray*}
\lim_{n \to \infty} \mathrm e^{-r \tau\wedge \tau_n} W(X_{\tau\wedge \tau_n},M_{\tau\wedge \tau_n})1_{\{\tau<\infty\}} = \mathrm e^{-r \tau} W(X_{\tau},M_\tau) 1_{\{\tau< \infty \}},
\end{eqnarray*}
$\mathbf P_{x,m}$-almost surely and, hence, in $\mathbf P_{x,m}$-probability, and because, by Corollary \ref{zoe}, the sequence $(\mathrm e^{-r \tau\wedge \tau_n} W(X_{\tau\wedge \tau_n},M_{\tau\wedge \tau_n})1_{\{\tau<\infty\}})_{n \in \mathbb N}$ is uniformly integrable, we have
\begin{eqnarray}
\lim_{n \to \infty} \mathbf E _{x,m}\hskip 0.3mm [\mathrm e^{-r \tau\wedge \tau_n}W(X_{\tau\wedge \tau_n},M_{\tau\wedge \tau_n})1_{\{\tau<\infty\}}] = \mathbf E_{x,m} \hskip 0.3mm [ \mathrm e^{-r \tau} W(X_{\tau},M_\tau)1_{\{\tau<\infty\}}] \label{etde2}
\end{eqnarray}
by Vitali's convergence theorem. Finally, over $\{\tau=\infty\}$, we have $\mathrm e^{-r \tau\wedge \tau_n} W(X_{\tau\wedge \tau_n},M_{\tau\wedge \tau_n}) = \mathrm e^{-r \tau_n} W(X_{\tau_n},M_{\tau_n})$. For $n$ large enough, $x \in (\alpha_n, \beta_n)$. Therefore,
\begin{align*}
 \mathbf E_{x,m}\hskip 0.3mm[ \mathrm e^{-r \tau_n} & W(X_{\tau_n},M_{\tau_n})]
\\
&= \mathbf E_{x,m} \hskip 0.3mm [\mathrm e^{-r\tau_n} W(X_{\tau_n},M_{\tau_n}) 1_{\{X_{\tau_n}  = \alpha_n\}}]+ \mathbf E_{x,m}\hskip 0.3mm [\mathrm e^{-r \tau_n} W(X_{\tau_n} ,M_{ \tau_n})1_{\{X_{\tau_n}  = \beta_n\}}] 
\\
&\leq  \mathbf E _{x,m}\hskip 0.3mm [\mathrm e^{-r\tau(\alpha_n)} W(\alpha_n, M_{\tau(\alpha_n)})]+ \mathbf E_{x,m}\hskip 0.3mm [\mathrm e^{-r \tau(\beta_n)}W( \beta_n,\beta_n)] \allowdisplaybreaks
\\
&\leq [1-F(m)] \! \left[ {h_2(x) \over h_2(\alpha_n)} \, G(\alpha_n) + {h_1(x) \over h_1(\beta_n)}\, G(\beta_n) \right]
\end{align*}
for any such $n$, where the second inequality follows from Lemma \ref{lem_uniform_integrability} along with the fact that the maximum process $M$ is nondecreasing. Together with Corollary \ref{corosigne} and the growth property (\ref{gp'}), this implies
\begin{eqnarray}
\lim_{n \to \infty}\mathbf E_{x,m}\hskip 0.3mm[\mathrm e^{-r \tau_n} W(X_{\tau_n},M_{\tau_n})] =0. \label{etde3}
\end{eqnarray}
Using (\ref{etde1})--(\ref{etde3}) to take the limit as $n$ goes to $\infty$ in \eqref{eq_ito_dynkin_1} yields (\ref{superharmonic}). The result follows. \hfill $\blacksquare$

\subsection{Wrapping Up} \label{Sreemati}

We are now ready to complete the proof of Proposition \ref{=} and, thereby, of Theorem \ref{maintheo}. First, by (\ref{mainproblem}) and \eqref{defW}, we have $V_b \leq V$. Next, by Lemmas \ref{lem_above} and \ref{lem_superharmonic}, we have
\begin{eqnarray*}
W(x,m)\geq \mathbf E_{x,m} \! \left[[1- F(M_ \tau)]\, \mathrm e^{-r \tau} R(X_\tau)  +\int_{0}^\tau \mathrm e^{-r t}G(M_t)f(M_t) \, \mathrm dM_t \right]\hskip -1mm,
\end{eqnarray*}
for all $\tau \in\mathcal T_{x,m}$ and $(x,m) \in \mathcal J$. Taking the supremum over $\tau$, we obtain $V \leq W$ by \eqref{mainproblem}. All we then need is the following lemma.

\begin{lem}
$W=V_b$.
\end{lem}

\noindent{\textbf{Proof.} By (\ref{Z1}) and (\ref{defW}), we have $W (x,m)= V_b(x,m)$ for all $(x,m) \in \mathcal S$. Now, let $(x,m) \in \mathcal C$ and, for each $n \in \mathbb N$, define $\tau_n$ as in the proof of Lemma \ref{lem_superharmonic}. For each $t < \tau_b \wedge \tau_n$, we have $(X_t, M_t) \in \mathcal C$, $\mathbf P_{x,m}$-almost surely, and thus $X_t \notin \{ x_R, b(M_t)\}$. Hence, applying \eqref{itoforW} to $\tau_f = \tau_b \wedge \tau_n$, we have
\begin{align*}
W(x,m)  = \;& \mathrm e^{-r\tau_b\wedge \tau_n} W(X_{\tau_b\wedge \tau_n},M_{\tau_b\wedge \tau_n}) + \int_0^{\tau_b\wedge \tau_n} \mathrm e^{-rt} (r  W- \mathcal{L}  W) (X_t,M_t) \, \mathrm dt \allowdisplaybreaks
\\
&\!\! + \int_0^{\tau_b\wedge \tau_n} \mathrm e^{-rt}G(M_t)f(M_t) \, \mathrm dM_t - \int_0^{\tau_b\wedge \tau_n} \mathrm e^{-rt}\, \frac{\partial W}{\partial x}\,(X_t,M_t) \, \sigma(X_t) \, \mathrm dW_t,
\end{align*}
$\mathbf P_{x,m}$-almost surely. The second term on the right-hand side vanishes as $\mathcal{L} W -rW=0$ over $\mathcal C$, and the fourth term on the right-hand side has zero expectation as shown in the proof of Lemma \ref{lem_superharmonic}. Therefore,
\begin{eqnarray*}
W(x,m)= \mathbf E_{x,m}\!\left[\mathrm e^{-r\tau_b\wedge \tau_n} W(X_{\tau_b\wedge \tau_n},M_{\tau_b\wedge \tau_n}) + \int_0^{\tau_b\wedge \tau_n} \mathrm e^{-rt}G(M_t)f(M_t) \, \mathrm dM_t \right]\hskip -1mm.
\end{eqnarray*}
Letting $n$ go to $\infty$ as in the proof of Lemma \ref{lem_superharmonic}, we obtain
\begin{align*}
W(x,m) &=  \mathbf E_{x,m}\!\left[\mathrm e^{-r\tau_b} W(X_{\tau_b},M_{\tau_b}) + \int_0^{\tau_b} \mathrm e^{-rt}G(M_t)f(M_t) \, \mathrm dM_t \right]
\\
&= \mathbf E_{x,m}\!\left[[1- F(M_{\tau_b})]\,\mathrm e^{-r\tau_b} R(X_{\tau_b}) + \int_0^{\tau_b} \mathrm e^{-rt}G(M_t)f(M_t) \, \mathrm dM_t \right]
\\
&=V_b(x,m)
\end{align*}
by (\ref{Z1}) and \eqref{defW}. The result follows. \hfill $\blacksquare$

\bigskip

The proof of Theorem \ref{maintheo} is now complete.

\section{Discussion} \label{Discussion}

In this section, we discuss the implications of our model in the context of the motivating example introduced in Section \ref{AMotExa}.

\subsection{Properties of the Value Function}

Intuitively, an increase in $m$, holding $x$ fixed, brings bad news for the DM unless it is accompanied by a breakthrough, for it means that he will have to wait for $X$ to reach a higher maximum value before hoping to benefit from a breakthrough. Our next result confirms this intuition, and moreover shows that the marginal cost of an increase in $m$ is discontinuous across the open half-line $(\alpha, \underline m) \times \{\underline m\}$.

\begin{pro} \label{propertyV1}
The following holds$:$
\begin{itemize}

\item[(i)]

${\partial V\over \partial m}<0$ over $\mathcal J \setminus (\mathcal J_{\underline m} \setminus \{(\underline m, \underline m)\})$.

\item[(ii)]

For each $(x, \underline m) \in \mathcal J_{\underline m} \setminus \{(\underline m, \underline m)\},$ ${\partial V \over \partial m}\, (x, \underline m^+) > {\partial V \over \partial m}\, (x, \underline m^-)$.
\end{itemize}
\end{pro}

Hence, over the half-line $(\alpha,\underline m] \times \{\underline m\}$, $\partial V \over \partial m$ is continuous at the endpoint $(\underline m, \underline m)$---because the Neumann condition holds everywhere over $\mathcal D$---but not at the points $(x,\underline m)$ for $x < \underline m$. This notably reflects that, whereas $B(\underline m)= 0$---which is required by the continuous-fit condition---we have $B'(\underline m) >0$. One can also check that the Neumann condition at $(\underline m, \underline m)$ implies $A'(\underline m) < C'(\underline m)$. This discontinuity intuitively reflects the change of regime across the half-line $(\alpha,\underline m] \times \{\underline m\}$: before $X$ reaches $\underline m$, the DM only invests in the superior technology, should it become available; afterwards, the DM starts investing in the stand-alone technology, which happens if $X$ decreases enough after having reached a maximum.

We now examine the impact on the DM's value of an increase in $x$, holding $m$ fixed. To obtain a clear-cut result, we have to make additional assumptions about the diffusion $X$ and the payoff function $R$. The following result holds.

\begin{pro} \label{propertyV2}
If $h_1$ and $h_2$ are convex and $R'>0$ over $\mathcal I,$ then ${\partial V\over \partial x}>0$ over $\mathcal J$.
\end{pro}

The assumption that $h_1$ and $h_2$ be convex is mild and is typically satisfied in real-options models of investment; for instance, this is the case if $X$ follows a geometric Brownian motion. More generally, because the lower endpoint $\alpha$ of $\mathcal I$ is, by assumption, inaccessible, a sufficient condition is that the net depreciation $rx - \mu(x)$ of an asset yielding a revenue flow $X$ be nondecreasing in $x$ (Alvarez (2003, Corollary 1)).

\subsection{Implications for Investment Theory}

A key feature of the optimal investment strategy is that investment in the stand-alone technology with payoff function $R$ only takes place when the process $(X,M)$ hits the free boundary $x= b(m)$, that is, when, after having reached a maximum value $m$, the process $X$ drops down to the lower threshold $b(m)$. Thus, unlike in the standard real-options model (Dixit and Pindyck (1994)), investment in the stand-alone technology takes place in busts rather than in booms; this reflects that the stand-alone technology becomes attractive for the DM only if he becomes sufficiently pessimistic about eventually benefiting from the superior technology with payoff function $U$. In line with this intuition, a necessary condition for investment in the stand-alone technology is that $X$ must have reached $\underline m$.

By contrast, if $U$ satisfies A1--A3, then investment in the superior technology only takes place, after a breakthrough has occurred, when $X$ reaches an upper threshold $x_U$; that is, in booms rather than in busts. For instance, in the specification outlined in Section \ref{TecAss} and detailed in Appendix C, we have $x_U > x_R$; hence, when $m \geq \underline m$, either $X$ drops down to $b(m)$ before reaching a new maximum value, which triggers investment in the stand-alone technology, or, if a breakthrough occurs when $X$ reaches a new maximum value, investment in the superior technology takes place immediately.

This qualitative difference between investments taking place in booms and investments taking place in busts is a novel prediction of the model. An empirical implication is that investments requiring cooperation from outside developers should take place in booms, in contrast with investments involving a more routine technology. A further implication of our model is that, because $b(m) >x_R$ is strictly increasing in $m > \underline m$, the return required to invest in the stand-alone technology is path-dependent and is always higher than it would be, were this the only available technology. In particular, when $\beta = \infty$ and $R$ is strictly increasing and unbounded above, the limit condition \eqref{tc} implies that the required return for investing in the stand-alone technology may assume arbitrarily high values. Intuitively, this is because, under technological uncertainty, the DM is willing to further delay investment in the stand-alone technology, and requires a higher return to give up the option to invest in the superior technology.

It is interesting to compare this pattern of investment to that arising in models of investment under uncertainty where the DM can invest in alternative projects. Building on Dixit (1993), D\'ecamps, Mariotti, and Villeneuve (2006) suppose that the DM can invest in two alternative projects, one with a high investment cost and a high output rate, and the other with a low investment cost and a low output rate. They show that investment in the former only takes place in booms, while, for certain values of the parameters, investment in the latter can also take place in busts, that is, when the DM becomes sufficiently pessimistic about the output price recovering soon enough to make investment in the high-cost and high-output project worthwhile again. The distinguishing feature of the present model is that the superior technology is not present from the outset, as it is only supplied by the developers when the share of the surplus they can secure by bargaining with the DM covers their development cost; moreover, the time at which such a breakthrough occurs is unknown to the DM because the development cost is the developers' private information. Thus the evolution of the cash-flow process $X$ provides information both about the desirability of investment and about the developers' cost, leading to a rich two-dimensional dynamics. A similar duality arises in the war-of-attrition investment model of D\'ecamps and Mariotti (2004), where each player does not observe his rival's cost; however, their model is cast in a Poisson rather than in a Brownian framework, which allows one to reduce each player's decision problem to a standard one-dimensional optimal stopping problem.

In our model, investment in the stand-alone technology optimally takes place in busts because the DM faces a downside risk as well as a standard upside potential. This feature also arises when the DM learns about the drift of the cash-flow process, as in D\'ecamps, Mariotti, and Villeneuve (2005) and Klein (2009). The difference is that the downside risk in our model is not tied to an intrinsic but unknown characteristic of the investment project, but rather to the fact that, while the DM has the option to wait until a superior technology becomes available, he may have to resign himself to use the stand-alone technology should the underlying cash-flow deteriorate too much.

\subsection{Comparative Statics}

Our comparative-statics results rely on two partial orders over the sets of development-cost distributions and of payoff functions that are motivated by (\ref{E}) and (\ref{edocase1a})--(\ref{edocase1ain}). We throughout assume that $U$ satisfies the same assumptions as $R$, so that, following a breakthrough, it is optimal for the DM to invest as soon as $X$ reaches a threshold $x_U$, and that $P' >0$ over $\mathcal I$. A case in point is the specification outlined in Section \ref{TecAss}.

First, let $F_{Z,1}$ and $F_{Z,2}$ be two distributions of development costs with $\mathcal C^1$ densities $f_{Z,1}>0$ and $f_{Z,2}>0$ over $\mathbb R_+$. Following Shaked and Shanthikumar (2007, Section 1.B.1), we say that $F_{Z,2}$ \textit{dominates} $F_{Z,1}$ \textit{in the hazard-rate order} if
\begin{eqnarray}
{f_{Z,1} \over 1- F_{Z,1}}\, > {f_{Z,2}  \over 1- F_{Z,2}} \label{distup}
\end{eqnarray}
over $\mathbb R_+$, so that, in particular, $F_{Z,2}$ first-order stochastically dominates $F_{Z,1}$: development costs tend to to be higher under $F_{Z,2}$ than under $F_{Z,1}$. The following result then holds.

\begin{pro} \label{compstat1}
Let $b_1: [\underline m \,\!_1, \beta) \to [x_R, \beta)$ and $b_2: [\underline m \,\!_2, \beta) \to [x_R, \beta)$ be the optimal free boundaries under the distributions $F_{Z,1}$ and $F_{Z,2}$ of development costs. If $F_{Z,2}$ dominates $F_{Z,1}$ in the hazard-rate order$,$ then $\underline m \,\!_1 > \underline m \,\!_2$ and $b_2 > b_1$ over $[\underline m \,\!_1, \beta)$.
\end{pro}

A consequence of this result is that, if the development cost increases in the hazard-rate order, then the DM becomes more cautious and less prone to bear downside risk. As a result, he is more likely to invest in the stand-alone technology.

Next, let $U_1$ and $U_2$ two payoff functions from investing in the new technology, and let $F_Z$ be the distribution of development costs. We say that $U_2$ \textit{dominates} $U_1$ if $U_2 > U_1$ and $U'_2 > U'_1$. As usual, we say that $F_Z$ satisfies the \textit{monotone hazard-rate property} (MHRP) if $f_Z \over 1- F_Z$ is nondecreasing over $\mathcal I$. The following result then holds.

\begin{pro} \label{compstat2}
Suppose that $h_1$ is convex$,$ $U_1$ is concave$,$ and $F_Z$ satisfies MHRP$,$ and let $b_1: [\underline m \,\!_1, \beta) \to [x_R, \beta)$ and $b_2: [\underline m \,\!_2, \beta) \to [x_R, \beta)$ be the optimal free boundaries under the payoffs functions $U_1$ and $U_2$ from investing in the new technology. If $U_2$ dominates $U_1,$  then $\underline m \,\!_2 > \underline m \,\!_1$ and $b_1 > b_2$ over $[\underline m \,\!_2, \beta)$.
\end{pro}

When the payoff function from investing in the new technology increases from $U_1$ to $U_2$ in the partial order we have defined, two effects are at play. The direct effect is that, by (\ref{at}), the value $V_U$ increases, and hence that so does, by (\ref{NASH}), the continuation value $G$ that the DM can obtain by bargaining with a successful developer. The indirect effect operates through the incentives of the developers. Notice, indeed, that the breakthrough rate $H$ in (\ref{E}) can be written as $H = P' \, {f_Z \circ P \over 1- F_Z \circ P}$, and hence depends on the share $P$ of the surplus $V_U - V_R$ that a successful developer can obtain by bargaining with the DM, as well as on its derivative $P'$. If $h_1$ is convex and $U_1$ is concave, then, when $U$ increases from $U_1$ to $U_2$, both $P$ and $P'$ increase, reflecting that the developers stand to gain both in absolute and marginal terms;\footnote{The convexity assumptions on $h_1$ and $U_1$ are not needed for this result if the optimal investment threshold is higher under $U_2$ than under $U_1$.} if, moreover, $F_Z$ satisfies MHRP, then we obtain that the breakthrough rate increases, as in Proposition \ref{compstat1}. Overall, an increase in $U$ makes it more likely that the DM will benefit from an even superior technology earlier on, which makes him more prone to bear downside risk; observe that the indirect effect described above is absent from models that treat technological breakthroughs as exogenous.

\section{Concluding Remarks} \label{Concluding Remarks}

In this paper, we have provided a new model of investment under technological and cash-flow uncertainty. The distinctive feature of our model is that the values of the stand-alone technology and of the superior technology depend on current market conditions, and that the occurrence of technological breakthroughs is correlated with the evolution of market conditions; hence, in our motivating example, a new technology is introduced as soon as market conditions are favorable enough to make it profitable for developers to do so, given the share of the surplus from the innovation that accrues to a successful developer. Thus our model may be seen as a first step towards a theory of the interactions between both sides of the market for technological innovations.

The main insight from our analysis is that investment in the stand-alone technology should only occurs in busts, when the market conditions deteriorate enough after having reached a maximum; we provide a complete characterization of the corresponding optimal investment boundary. By contrast, investments in new technologies requiring the active cooperation of developers should take place in booms. This intuitively reflects that the stand-alone technology becomes attractive only when the firm becomes pessimistic enough about a breakthrough shortly being forthcoming. As a result, and in contrast with standard models of investment under uncertainty, the firm bears downside risk, in addition to the upside potential associated to technological breakthroughs. A decrease in development costs, or an increase in the value of the new technology, makes the firm more prone to bear such downside risk and to delay investment in the stand-alone technology.

We have throughout assumed that at most two technologies are available. In that respect, our analysis is less rich than that of Balcer and Lippman (1984), who consider a sequence of innovations arising according to a semi-Markov process. We conjecture, however, that qualitatively similar results would hold if the model were extended to allow for multiple innovations. Another limitation of our analysis is that we have considered the investment policy of an isolated firm. A fascinating avenue of research would be to investigate the implication of technological and cash-flow uncertainty for the equilibrium of an industry. We leave these questions for future work.

\newpage

\small \baselineskip =5mm

\renewcommand{\thesection}{Appendix A: Ommitted Proofs}

\section{}

\newtheorem{applemma}{Lemma}[section]

\renewcommand{\thesection}{A}

\renewcommand{\theequation}{A.\arabic{equation}}
\setcounter{equation}{0}

\noindent{\textbf{Proof of Lemma \ref{lem_sign}.} By (\ref{sa}) and (\ref{laplace}),
\begin{eqnarray*}
V_R(x) \geq \mathbf E_x \hskip 0.3mm [\mathrm e^{-r \tau(y)} R(X_{\tau(y)}) ] = \frac{h_1(x)}{h_1(y)}\,R(y)
\end{eqnarray*}
for all $x \in \mathcal I$ and $y \in [x, \beta)$. Letting $y$ go to $\beta^-$ and taking advantage of (\ref{gp}), we obtain that $V_R \geq 0$ over $\mathcal I$ and, as $ R = V_R $ over $[x_R, \beta)$ by (\ref{VR}), that $R \geq 0$ over $[x_R, \beta)$. To show that these inequalities are strict, observe from A3 that $R$ cannot be identically zero over $[x_R, \beta)$. Thus $R(y) >0$ for some $y \in [x_R, \beta)$. Because, by (\ref{sa}), (\ref{laplace}), and (\ref{VR}),
\begin{eqnarray*}
R(x) = V_R(x) \geq  \mathbf E_x\hskip 0.3mm [\mathrm e^{-r\tau(y)} R(X_{\tau(y)}) ]  =  \left\{ \begin {matrix} \frac{h_1(x)}{h_1(y)}\,R(y) & \text{if} & x \leq y, \\  \frac{h_2(x)}{h_2(y)} \, R(y) & \text{if} & x > y, \end{matrix} \right.
\end{eqnarray*}
for all $x \in [x_R, \beta)$, this implies that $R >0$ over $[x_R, \beta)$, which, along with (\ref{VR}) again, in turn implies that $V_R >0$ over $\mathcal I$. The result follows. \hfill $\blacksquare$

\bigskip

\noindent{\textbf{Proof of Lemma \ref{explodediadlem}.} For each $m \in [x_R, \beta)$, we have $G(m)>R(m)$ by A6, and $L(m)<0$ as shown in the text. It follows that
\begin{eqnarray*}
\lim_{(x,m'), x<m' \rightarrow (m,m)} D(x,m')E(x,m') = \frac{H(m) \sigma^2(m)}{2 L(m)h_2(m)} \, \gamma S'(m)h_2(m)[R(m)-G(m)] >0.
\end{eqnarray*}
Thus the mapping $(x,m) \mapsto D(x,m)E(x,m)$ can be continuously extended over $\overline{\mathcal J_E}$, which implies (\ref{explodediad}) as $D$ vanishes over $\mathcal D$. The result follows. \hfill $\blacksquare$

\bigskip

\noindent{\textbf{Proof of Lemma \ref{odeforb}.} To avoid notational clutter, we hereafter write $b$ for $b(m)$ and $b'$ for $b'(m)$. Notice first that, over the subset $\{(x,m) \in \mathcal J: m \geq \underline m \mbox{ and } x \in (b,m]\}$ of $\mathcal C$, any solution $W \in \mathcal V$ to (\ref{vs1}) is of the form
\begin{eqnarray}
W(x,m) = A(m) h_1(x) + B(m) h_2 (x) \label{GW}
\end{eqnarray}
for some functions $A,B \in \mathcal C^1([\underline m, \beta))$. By (\ref{GW}), the continuous- and smooth-fit conditions $W(b,m) = [1-F(m)]R(b)$ and $\frac{\partial W}{\partial x}\,(b, m) =[1-F(m)]R'(b)$ are satisfied if and only if
\begin{align}
A(m)h_1(b)+B(m)h_2(b)  &= [1-F(m)]R(b), \label{hahaha}
\\
A(m)h'_1(b)+B(m)h'_2(b)  &= [1-F(m)]R'(b). \label{hihihi}
\end{align}
By (\ref{hihihi}), we have
\begin{eqnarray*}
B(m)  =  \frac{1}{h'_2(b)}\,\{[1-F(m)]R'(b) - A(m)h'_1(b)\}.
\end{eqnarray*}
Substituting in (\ref{hahaha}) and multiplying by $h'_2(b)$, we obtain
\begin{eqnarray*}
A(m)h_1(b)h'_2(b)+[1-F(m)]R'(b)h_2(b) - A(m)h'_1(b)h_2(b) = [1-F(m)]R(b)h'_2(b)
\end{eqnarray*}
and thus, by (\ref{wronskian}),
\begin{eqnarray}\label{GA}
A(m) = \frac{1-F(m)}{\gamma S'(b)} \,[R'(b)h_2( b)- R(b) h'_2(b)].
\end{eqnarray}
By symmetry,
\begin{eqnarray}\label{GB}
B(m) = - \frac{1-F(m)}{\gamma S'(b)} \,[R'(b)h_1( b)- R(b) h'_1(b)].
\end{eqnarray}
By (\ref{GW}), the Neumann condition $\frac{\partial W}{\partial m}\,(m,m) =-f(m)G(m)$ is satisfied if and only if
\begin{eqnarray}
A'(m)h_1(m)+B'(m)h_2(m) = -f(m)G(m). \label{hohoho}
\end{eqnarray}
Differentiating (\ref{GA}) with respect to $m$ yields
\begin{align}
A'(m) = \;& \frac{b' [1-F(m)]}{\gamma S'(b)}\, \bigg\{\underbrace{R''(b)h_2(b)- R(b)h_2''(b)- \frac{S''(b)}{S'(b)}\,[ R'(b)h_2( b)- R(b) h'_2(b) ]  }_{\large{ \equiv Q(b)}} \bigg\} \notag
\\
&\! -\frac{f(m)}{\gamma S'(b)}\,[ R'(b)h_2( b)- R(b) h'_2(b)].  \label{GA'''}
\end{align}
Now, observe that
\begin{align}
Q(b) &=   \frac{2}{\sigma^2(b)} \left\{\frac{1}{2}\, \sigma^2(b) R''(b)h_2(b)- \frac{1}{2} \, \sigma^2(b) R(b)h_2''(b) + \mu(b)[ R'(b)h_2( b)- R(b) h'_2(b) ] \right\} \notag \\
&= \frac{2}{\sigma^2(b)}\,[(\CL R -r R)(b)h_2(b) - R(b)(\CL h_2 - rh_2)(b)] \notag
\\
&=  \frac{2}{\sigma^2(b)}\, (\CL R -r R)(b) h_2(b) \notag
\\
&=  \frac{2}{\sigma^2(b)}\, L(b) h_2(b), \label{forneumann}
\end{align}
where the first equality follows from noticing that $\frac{S''(b)}{S'(b)}= - \frac{2 \mu(b)}{\sigma^2(b)}$ by \eqref{scale}, the second equality follows from (\ref{defmathcalL}), the third equality follows from the fact that $\CL h_2 - rh_2 = 0$, and the fourth equality follows from \eqref{L}. Using (\ref{forneumann}) to rewrite (\ref{GA'''}) yields
\begin{eqnarray}
A'(m) = \frac{1}{\gamma S'(b)} \left\{\frac{2b' [1-F(m)]}{\sigma^2(b)}\, L(b) h_2(b)- f(m)[ R'(b)h_2( b)- R(b) h'_2(b)] \right\}\!.  \label{GA''''}
\end{eqnarray}
By symmetry,
\begin{eqnarray}
B'(m) = - \frac{1}{\gamma S'(b)} \left\{\frac{2b' [1-F(m)]}{\sigma^2(b)}\, L(b) h_1(b)- f(m)[ R'(b)h_1( b)- R(b) h'_1(b)] \right\}\!.  \label{GB''''}
\end{eqnarray}
Using (\ref{D})--(\ref{H}), we obtain from (\ref{GA''''})--(\ref{GB''''}) that (\ref{hohoho}) holds if and only if
\begin{eqnarray*}
\frac{2b'}{\sigma^2(b)}\, L(b) D(b,m) = H(m) \{ R'(b)D(b,m) +R(b) [h'_1(b) h_2(m)- h_2'(b)h_1(m)]-\gamma S'(b) G(m)\},
\end{eqnarray*}
that is, observing that
\begin{align*}
h'_1(b) h_2(m)- h_2'(b)h_1(m) &=  {\gamma S'(b) +h_1(b)h'_2(b) \over h_2(b)} \, h_2(m)- h_2'(b)h_1(m)
\\
&= {\gamma S'(b)h_2(m) - h_2'(b)D(b,m) \over  h_2(b)}
\end{align*}
by (\ref{wronskian}) and (\ref{D}), and that $D(b,m) >0$ as $b <m$, if and only if
\begin{eqnarray*}
b' = \frac{H(m)\sigma^2(b)}{2L(b)h_2(b) }  \left\{  \frac{\gamma S'(b)}{D(b,m)}\,[R(b)h_2(m)-G(m)h_2(b)] + R'(b)h_2( b)- R(b) h'_2(b) \right\}\!,
\end{eqnarray*}
which is \eqref{edocase1a} by \eqref{E}. The result follows. \hfill $\blacksquare$

\bigskip

\noindent \textbf{Proof of Lemma \ref{WinF}.} First, $W \in \mathcal C^0(\mathcal J)$ by construction. Next, as $R$, $h_1$, and $h_2$ are $\mathcal C^2$, and $F$, $G$, and $b$ are $\mathcal C^1$, the functions (\ref{Z1})--(\ref{Z3}) and (\ref{Z4}) are $\mathcal C^{2,1}$ over the domains $\mathcal S$, $\mathcal J_1=(\alpha, x_R] \times [\underline  m, \beta )$, $\mathcal J_2=\{(x,m) \in \mathcal J: m \geq \underline m \mbox{ and } x \in  [b(m), m]\}$ and $\mathcal J_3=\{(x,m) \in \mathcal J : m \leq \underline m \}$, respectively. Finally, that $W \in \mathcal C^1(\mathcal J \setminus \mathcal J_{\underline m})$ follows from the following observations. First, the function obtained by pasting together (\ref{Z1})--(\ref{Z2}) is $\mathcal C^1$ at $\{x_R\} \times (\underline m, \beta)$ by the smooth-fit property for $V_R$ at $x_R$; Second, the function obtained by pasting together (\ref{Z1}) and (\ref{Z3}) is $\mathcal C^1$ at $\{(x,m) \in \mathcal J: m > \underline m \mbox{ and } x=b(m)\}$ because the ``horizontal'' smooth-fit condition
\begin{eqnarray*}
A(m)h_1'(b(m)) + B(m) h_2'(b(m)) = {\partial W \over \partial x} \, (b(m)^+\hskip -0.5mm,m) = {\partial W \over \partial x} \, (b(m)^-\hskip -0.5mm,m) = [1- F(m)] R'(b(m))
\end{eqnarray*}
implies, upon differentiating the continuous-fit condition
\begin{eqnarray*}
W(b(m),m)= A(m)h_1(b(m)) + B(m) h_2(b(m)) = [1- F(m)] R(b(m)),
\end{eqnarray*}
the ``vertical'' smooth-fit condition
\begin{eqnarray*}
{\partial W \over \partial m} \, (b(m),m^-) =A'(m)h_1(b(m)) + B'(m)h_2(b(m)) = -f(m) R(b(m)) ={\partial W \over \partial m} \, (b(m),m^+).
\end{eqnarray*}
Thus $W \in \mathcal V$, as claimed. Notice finally that the function obtained by pasting together (\ref{Z3}) and (\ref{Z4}) has continuous partial derivatives at $(\underline m, \underline m)$ by the Neumann condition along with the fact that $A(\underline m)= C(\underline m)$. The result follows. \hfill $\blacksquare$

\bigskip

\noindent \textbf{Proof of Corollary \ref{corosigne}.} Recall from A6 and Lemma \ref{lem_sign} that $R >0$ over $[x_R, \beta)$ and that $G >0$ over $\mathcal I$. We consider four cases in turn.

\subparagraph{Case 1}

Suppose first that $(x,m) \in \mathcal S$, so that $W(x,m)$ is given by (\ref{Z1}). Then $x \geq x_R$ and thus $W(x,m) >0$ as $R >0$ over $[x_R, \beta)$.

\subparagraph{Case 2}

Suppose next that $(x,m) \in (\alpha,x_R) \times [\underline m, \beta)$, so that $W(x,m)$ is given by (\ref{Z2}). Then $W(x,m) >0$ as $R (x_R)>0$.

\subparagraph{Case 3}

Suppose now that $m \in [\underline m, \beta)$ and $x \in  (b(m), m]$, so that $W(x,m)$ is given by (\ref{Z3})--(\ref{exprB}). We just need to check that $A(m)$ and $B(m)$ are nonnegative, with one of them strictly positive. By (\ref{23:42}) and \eqref{exprA}, we have
\begin{eqnarray}
A(m) = [1-F(m)] \hat R' (\zeta(b(m)) >0 \label{brownie}
\end{eqnarray}
because $\zeta(b(m)) \geq \zeta(x_R)$ and $\hat R' >0$ over $[\zeta(x_R), \infty)$ as shown in the proof of Lemma \ref{mx}. By (\ref{23:42}) and \eqref{exprB}, we have $B(\underline m)=0$, and for $m \in(\underline m,\beta)$, we have
\begin{align*}
B(m) &=  [1-F(m)] \zeta(b(m))\! \left[{R(b(m)) \over h_1(b(m))} - \hat R'(\zeta(b(m)))\right]
\\
&>  [1-F(m)] \zeta(b(m))\! \left[{R(b(m)) \over h_1(b(m))} - {\hat R(\zeta(b(m))) - \hat R(\zeta(x_R)) \over \zeta(b(m))- \zeta(x_R)}\right]
\\
&=  [1 - F(m)] \,{ \zeta(b(m)) \zeta(x_R) \over  \zeta(b(m))- \zeta(x_R)}\!\left[ {R(x_R) \over h_1(x_R)} - {R(b(m)) \over h_1(b(m))}  \right]
\\
&\geq  0,
\end{align*}
where the first inequality follows from the strict concavity of $\hat R$ over $[\zeta(x_R), \infty)$, and the second inequality follows from the optimality of the stopping threshold $x_R$ for \eqref{sa}.

\subparagraph{Case 4}

Suppose finally that $m \in (\alpha, \underline m)$ and $x \in  (\alpha, m]$, so that $W(x,m)$ is given by (\ref{Z4})--(\ref{exprAbelow}). Then $W(x,m) >0$ as $R(x_R) > 0$ and $G > 0$ over $\mathcal I$. Hence the result. \hfill $\blacksquare$

\bigskip

\noindent \textbf{Proof of Proposition \ref{propertyV1}.} (i) That ${\partial V \over \partial m} \,(\underline m, \underline m)<0$ directly follows from the Neumann condition. We then consider four cases in turn.

\subparagraph{Case 1}

Suppose first that $(x,m) \in \mathcal S \setminus \{(x_R, \underline m)\}$, so that $V(x,m)$ is given by (\ref{Z1}). Then ${\partial V \over \partial m} \, (x,m) = -f(m) R(x) <0$ as $f>0$ over $\mathcal I$ and $R >0$ over $[x_R, \beta)$.

\subparagraph{Case 2}

Suppose next that $(x,m) \in (\alpha,x_R) \times (\underline m, \beta)$, so that $V(x,m)$ is given by (\ref{Z2}). Then ${\partial V \over \partial m} \, (x,m) = -f(m) \, {h_1(x) \over h_1(x_R) } \, R(x_R) <0$ as $f>0$ and $h_1>0$ over $\mathcal I$ and $R (x_R)>0$.

\subparagraph{Case 3}

Suppose now that $m \in (\underline m, \beta)$ and $x \in  (b(m), m]$, so that $V(x,m)$ is given by (\ref{Z3})--(\ref{exprB}). Then ${\partial V \over \partial m} \, (x,m) = A'(m) h_1(x) + B'(m) h_2(x)$. By \eqref{brownie}, we have
\begin{eqnarray}
A'(m) = -f(m) \hat R' (\zeta(b(m)) + [1-F(m)] \hat R'' (\zeta(b(m)) \zeta'(b(m)) b'(m)<0  \label{brownie'}
\end{eqnarray}
as  $f>0$ and $\zeta' >0$ over $\mathcal I$, $b'>0$ over $ [\underline m, \beta)$, and $\zeta(b(m)) \geq \zeta(x_R)$ and $\hat R' >0$ and $\hat R ''<0$ over $[\zeta(x_R), \infty)$ as shown in the proof of Lemma \ref{mx}. Because the mapping $x \mapsto {h_1( x) \over h_2(x)}$ is strictly increasing, \eqref{brownie'} implies that
\begin{align*}
{\partial V \over \partial m} \, (x,m)  &=  h_2(x) \!\left[ A'(m) \,{h_1(x) \over h_2(x)} + B'(m)\right]
\\
&\leq  h_2(x) \!\left[ A'(m) \,{h_1(b(m)) \over h_2(b(m))} + B'(m)\right]
\\
&= - {h_2(x) \over h_2(b(m))} \, f(m) R(b(m))
\\
&<  0,
\end{align*}
where the second equality follows from the vertical smooth-fit condition ${\partial V \over \partial m}\,(b(m), m) = -f(m) \linebreak R(b(m))$, and the second inequality follows from the fact that $f>0$ and $h_2 >0$ over $\mathcal I$ and that $R>0$ over $[x_R, \beta)$.

\subparagraph{Case 4}

Suppose finally that $m \in (\alpha, \underline m)$ and $x \in  (\alpha, m]$, so that $V(x,m)$ is given by (\ref{Z4})--(\ref{exprAbelow}). Then ${\partial V \over \partial m} \, (x,m) =-{f(m) \over h_1(m) } \, G(m)h_1(x) <0$ as $f>0 $, $h_1 >0$, and $G > 0$ over $\mathcal I$. This proves (i).

\vskip 3mm

(ii) We consider two cases in turn.

\subparagraph{Case 1}

Suppose first that $x \in (\alpha, x_R]$. On the one hand, by (\ref{Z1})--(\ref{Z2}), we have
\begin{eqnarray*}
{\partial V \over \partial m}\, (x, \underline m^+) = -f(\underline m) \,{h_1(x) \over h_1(x_R)} \, R(x_R).
\end{eqnarray*}
On the other hand, by (\ref{Z4})--(\ref{exprAbelow}), we have
\begin{eqnarray*}
{\partial V \over \partial m}\, (x, \underline m^-) = -f(\underline m) \,{h_1(x) \over h_1(\underline m)} \, G(\underline m),
\end{eqnarray*}
and the result follows from (\ref{comparpayoffs}).

\subparagraph{Case 2}

Suppose next that $x \in (x_R, \underline m)$. We first derive more compact expressions for $A'(m)$ and $B'(m)$, where $m \in(\underline m, \beta)$. To avoid notational clutter, we hereafter write $b$ for $b(m)$. By (\ref{Z3}), the vertical \pagebreak smooth-fit condition $\frac{\partial V}{\partial m}\,(b, m) =-f(m)R(b)$ and the Neumann condition $\frac{\partial V}{\partial m}\, (m,m) = -f(m)G(m)$ are satisfied if and only if
\begin{align}
A'(m)h_1(b) + B'(m)h_2(b)  &=  -f(m)R(b),  \label{haha}
\\
A'(m)h_1(m) + B'(m)h_2(m)  &= -f(m)G(m). \label{hihi}
\end{align}
By \eqref{hihi}, we have
\begin{eqnarray*}
B'(m)= - \frac{1}{h_2(m)}\,[f(m)G(m) + A'(m)h_1(m)].
\end{eqnarray*}
Substituting in (\ref{haha}) and multiplying by $h_2(m)$, we obtain
\begin{eqnarray*}
A'(m)h_1(b)h_2(m) -f(m)G(m)h_2(b) - A'(m)h_1(m)h_2(b)  =  -f(m)R(b)h_2(m)
\end{eqnarray*}
and thus, by (\ref{D}),
\begin{eqnarray}\label{GAA'}
A'(m)  =  \frac{f(m)}{D (b,m)}\, [R(b)h_2(m)  -G(m)h_2(b)].
\end{eqnarray}
By symmetry,
\begin{eqnarray}\label{GBB'}
B'(m)  = - \frac{f(m)}{D (b,m)}\, [R(b) h_1(m) - G(m)h_1(b) ].
\end{eqnarray}
Both (\ref{GAA'})--(\ref{GBB'}) can be extended by continuity at $\underline m$ to evaluate ${\partial V \over \partial m}\,(x,\underline m^+)$ for all $x \in (x_R, \underline m)$. Notice from (\ref{edocase1ain}) and (\ref{GBB'}) that
\begin{eqnarray*}
B'(\underline m)  = - \frac{f(\underline m)}{D (x_R,\underline m)}\, [R(x_R) h_1(\underline m) - G(\underline m)h_1(x_R) ] >0,
\end{eqnarray*}
taking again advantage from (\ref{comparpayoffs}). By (\ref{Z3}) and (\ref{Z4}), (i) holds at $x \in (x_R, \underline m)$ if and only if
\begin{eqnarray*}
A'(\underline m) h_1(x) + B'(\underline m) h_2(x) > C'(\underline m) h_1(x).
\end{eqnarray*}
Because $B'(\underline m) >0$ and the mapping $x\mapsto {h_2(x) \over h_1(x)}$ is strictly decreasing, this is the case for all $x \in (x_R, \underline m)$ if and only if
\begin{eqnarray}
A'(\underline m) h_1(\underline m) + B'(\underline m) h_2(\underline m) \geq C'(\underline m) h_1(\underline m). \label{jehaislesdimanches}
\end{eqnarray}
Using (\ref{D}), (\ref{edocase1ain}), (\ref{exprAbelow}), and (\ref{GAA'})--(\ref{GBB'}), it is easily checked that (\ref{jehaislesdimanches}) is in fact an equality. This, incidentally, reflects that $\partial V \over \partial m$ is continuous at $(\underline m, \underline m)$ because the Neumann condition holds everywhere over $\mathcal D$; notice also that this implies $A'(\underline m) < C'(\underline m)$ as $B'(\underline m)>0$. This proves (ii). Hence the result. \hfill $\blacksquare$

\bigskip

\noindent \textbf{Proof of Proposition \ref{propertyV2}.} We consider four cases in turn.

\subparagraph{Case 1}

Suppose first that $(x,m) \in \mathcal S$, so that $V(x,m)$ is given by (\ref{Z1}). Then ${\partial V \over \partial x} \, (x,m) = [1-F(m)] R'(x) >0$ as $R ' >0$ over $\mathcal I$.

\subparagraph{Case 2}

Suppose next that $(x,m) \in (\alpha,x_R) \times [\underline m, \beta)$, so that $V(x,m)$ is given by (\ref{Z2}). Then ${\partial V \over \partial x} \, (x,m) = [1-F(m)] \, {h_1'(x) \over h_1(x_R) } \, R(x_R) >0$ as $h_1' >0$ over $\mathcal I$ and $R (x_R)>0$.

\subparagraph{Case 3}

Suppose now that $m \in [\underline m, \beta)$ and $x \in  (b(m), m]$, so that $V(x,m)$ is given by (\ref{Z3})--(\ref{exprB}). Then ${\partial V \over \partial x} \, (x,m) = A(m) h_1'(x) + B(m) h_2'(x)$. It follows from the proof of Corollary \ref{corosigne} that $A(m) >0$ for $m \in[ \underline m, \beta)$ and that $B(m) >0$ for $m \in (\underline m, \beta)$, with $B( \underline m ) =0$. Hence $V(x,m)$ is convex in $x\in  (b(m), m]$ as $h_1$ and $h_2$ are convex. Thus ${\partial V \over \partial x} \, (x,m) >0$ because, by the smooth-fit condition, ${\partial V \over \partial x} \, (b(m),m)= [1-F(m)] R'(b(m))>0$ as $R' >0$ over $\mathcal I$.

\subparagraph{Case 4}

Suppose finally that $m \in (\alpha, \underline m)$ and $x \in  (\alpha, m]$, so that $V(x,m)$ is given by (\ref{Z4})--(\ref{exprAbelow}). Then ${\partial V \over \partial x} \, (x,m) =C(m)h_1'(x) >0$ as $h_1' >0$ over $\mathcal I$ and $C >0$ over $(\alpha, \underline m)$. Hence the result. \hfill $\blacksquare$

\bigskip

\noindent \textbf{Proof of Proposition \ref{compstat1}.} Consider the breakthrough densities $f_1 \equiv P'  f_{Z,1} \circ P $ and $f_2 \equiv P'  f_{Z,2} \circ P $ associated to $f_{Z,1}$ and $f_{Z,2}$. Then, as $P' >0$, the corresponding breakthrough rates $H_1$ and $H_2$ satisfy
\begin{eqnarray*}
H_1(m) \equiv P'(m ) \, {f_{Z,1}(P(m)) \over 1- F_{Z,1}(P(m))} > P'(m ) \,{f_{Z,2}(P(m)) \over 1- F_{Z,2}(P(m))} \equiv H_2(m)
\end{eqnarray*}
for all $m \in \mathcal I$. Using (\ref{E}) and Lemma \ref{mx}, we deduce from this that the vector fields $E_1$ and $E_2$ associated to $f_{Z,1}$ and $f_{Z,2}$ satisfy
\begin{eqnarray}
E_1(x,m) > E_2(x,m)  \label{elo}
\end{eqnarray}
for all $(x,m)$ such that $x \geq b_2(m)$. Now, suppose, by way of contradiction, that $b_1(m_0) \geq  b_2(m_0)$ for some $m_0 \in [\underline m \,\!_1, \beta) \cap [\underline m \, \! _2, \beta)$. Then, by (\ref{elo}), we have $b_1 (m) > b_2(m)$ for all $m \in ( m_0, \beta)$. Fix some $\varepsilon \in (0, \beta - m_0)$ and consider the following ODE:
\begin{align*}
b'(m)   &=  E_2(b (m),m), \quad m \geq m_0 + \varepsilon,
\\
b(m_0 + \varepsilon)  & =  b_1(m_0 + \varepsilon),
\end{align*}
with maximal solution $b_{12}$ satisfying $(b_{12}(m),m) \in \mathcal J_{E_2}^+$ for all $m$ in a maximal interval with lower endpoint $m_0 + \varepsilon$. By (\ref{elo}) again, we have $b_1 \geq b_{12} > b_2$ over this maximal interval, which must thus coincide with $[m_0 + \varepsilon, \beta)$. But then the interval $\mathcal I^0_2$ of possible endpoints for a candidate for the free boundary when the breakthrough density is $f_2$ cannot be reduced to a point, a contradiction. It follows that $b_2(m) > b_1(m)$ for all $m \in [\underline m \,\!_1, \beta) \cap [\underline m \, \! _2, \beta)$, and, in particular, that $\underline m \, \! _1 > \underline m\, \! _2$. Hence the result. \hfill $\blacksquare$

\bigskip

\noindent \textbf{Proof of Proposition \ref{compstat2}.} The proof consists of two steps.

\subparagraph{Step 1}

Consider the value functions $V_{U_1}$ and $V_{U_2}$ associated to $U_1$ and $U_2$ as in \eqref{saU}. Because $U_1$ and $U_2$ satisfy the same properties as $R$, we obtain
\begin{eqnarray*}
V_{U_1}(x) = \left\{ \begin{array}{lll} \frac{h_1(x)}{h_1(x_{U_1})}  \,U_1(x_{U_1})& \text{if} & x < x_{U_1}, \\  U_1(x) & \text{if} & x \geq x_{U_1}  \end{array} \right. \quad \mbox{and} \quad V_{U_2}(x) = \left\{ \begin{array}{lll} \frac{h_1(x)}{h_1(x_{U_2})}  \,U_2(x_{U_2})& \text{if} & x < x_{U_2}, \\  U_2(x) & \text{if} & x \geq x_{U_2}  \end{array}\right.
\end{eqnarray*}
for some optimal thresholds $x_{U_1}$ and $x_{U_2}$. We claim that $V_{U_2} > V_{U_1}$ and $V_{U_2}' > V_{U_1}'$ over $\mathcal I$. The first inequality follows directly from the assumption that $U_2 >U_1$. As for the second inequality, it is clear from the above expressions that it is satisfied over $(\alpha, x_{U_1} \wedge x_{U_2}]$ (the common part of the continuation regions) and over $[x_{U_1} \vee x_{U_2}, \beta)$ (the common part of the stopping regions). Consider now the region $(x_{U_1} \wedge x_{U_2}, x_{U_1} \vee x_{U_2})$. Suppose first that $x_{U_1} < x_{U_2}$. Then, for each $x \in (x_{U_1} ,x_{U_2})$,
\begin{align*}
V'_{U_2} (x)  - V'_{U_1}(x) &= {h'_1(x) \over h_1(x_{U_2})} \, U_2(x_{U_2})  - U'_1(x)
\\
&\geq  {h'_1(x_{U_1}) \over h_1(x_{U_2})} \, U_2(x_{U_2})  - U'_1(x_{U_1})
\\
&\geq   {h'_1(x_{U_1}) \over h_1(x_{U_1})} \, U_2(x_{U_1})  - U'_1(x_{U_1})
\\
& >  {h'_1(x_{U_1}) \over h_1(x_{U_1})} \, U_1(x_{U_1})  - U'_1(x_{U_1})
\\
&= 0,
\end{align*}
where the first inequality follows from the convexity of $h_1$ and the concavity of $U_1$, the second inequality follows from the optimality of the threshold $x_{U_2}$ under $U_2$, the third inequality follows from $U_2 > U_1$, and the second equality follows from the smooth-fit property under $U_1$.  Suppose next that $x_{U_2} < x_{U_1}$. Then, for each $x \in (x_{U_2} ,x_{U_1})$,
\begin{align*}
V'_{U_2} (x)  - V'_{U_1}(x) &=  U'_2(x)  - {h'_1(x) \over h_1(x_{U_1})} \,U_1(x_{U_1})
\\
&>  U'_1(x)  - {h'_1(x) \over h_1(x_{U_1})} \,U_1(x_{U_1}) \allowdisplaybreaks
\\
&\geq   U'_1(x_{U_1})  - {h'_1(x_{U_1}) \over h_1(x_{U_1})} \,U_1(x_{U_1})
\\
&=  0,
\end{align*}
where the first inequality follows  from $U'_2 > U'_1$, the second inequality follows from the convexity of $h_1$ and the concavity of $U_1$, and the second equality follows from the smooth-fit property under $U_1$. The claim follows.

\subparagraph{Step 2}

Together with \eqref{NASH}, Step 1 implies that $G_2 > G_1$, $P_2 > P_1$, and $P'_2 > P'_1$, with obvious notation. In particular, because $F_Z$ satisfies MHRP, the corresponding breakthrough rates $H_1$ and $H_2$ satisfy
\begin{eqnarray*}
H_2(m) = P_2'(m ) \, {f_Z(P_2(m)) \over 1- F_Z(P_2(m))} > P_1'(m ) \,{f_Z(P_1(m)) \over 1- F_Z(P_1(m))}= H_1(m)
\end{eqnarray*}
for all $m \in \mathcal I$. Using (\ref{E}) and Lemma \ref{mx}, we deduce from this along with $G_2 > G_1$ that the vector fields $E_1$ and $E_2$ associated to $U_1$ and $U_2$ satisfy \eqref{elo} for all $(x,m)$ such that $x \geq b_2(m)$. The remainder of the proof follows along the lines of the proof of Proposition \ref{compstat1}. Hence the result. \hfill $\blacksquare$

\renewcommand{\thesection}{Appendix B: The Dynamic Programming Principle}

\section{}

\renewcommand{\thesection}{B}

\renewcommand{\theequation}{B.\arabic{equation}}
\setcounter{equation}{0}

In this appendix, we show how to apply the dynamic programming principle to obtain the general form \eqref{mainpb-intro} of our problem, as announced in Section \ref{AMotExa}. On top of the assumptions made in Section \ref{TecAss}, we assume that $U \in \mathcal C^2(\mathcal I)$, with $U>R$ over $\mathcal I$, and that  $U$ satisfies A1--A3, so that there exists $x_U \in \mathcal I$ such that the stopping time $\tau_{X\geq x_U} \equiv \inf \hskip 0.5mm \{t \geq 0: X_t \geq x_U\}$ is the solution to the optimal stopping problem
\begin{eqnarray} \label{saU}
V_U(x) \equiv \sup_{\tau \in \CT_X} \mathbf E_x \hskip 0.3mm[ \mathrm e^{-r\tau } U(X_\tau)].
\end{eqnarray}
Recall that the payoff for the DM when stopping at $\tau \in \mathcal{T}_{X, X \geq Y}$  is
\begin{eqnarray}
\hat{J}(x,\tau) \equiv \overline{ \mathbf E}_x \! \left[1_{\{\tau<\tau_{X \geq Y}\! \}}\, \mathrm e^{-r \tau} R(X_{\tau})+1_{\{\tau \geq \tau_{X \geq Y} \! \}} \! \left[ \mathrm e^{-r \tau} U(X_{\tau} )  - \mathrm e^{-r \tau_{X \geq Y}} P(X_{\tau_{X \geq Y}}) \right]\right] \hskip -1mm.
\end{eqnarray}
The following result then holds.

\begin{lem}
For each $x\in \mathcal I$, $\sup_{\tau \in \mathcal{T}_{X , X \geq Y}} \hat {J} (x,\tau)= \overline V (x)$.
\end{lem}

\noindent \textbf{Proof.}
Given that, by (\ref{NASH}),
\begin{eqnarray*}
V_U(X_{\tau_{X\geq Y}}) - P(X_{\tau_{X \geq Y}}) = G(X_{\tau_{X\geq Y}}),
\end{eqnarray*}
it is sufficient to prove that
\begin{eqnarray*}
\sup_{\tau \in \mathcal{T}_{X, X \geq Y}} \hat{J}(x,\tau) =\sup_{\tau \in \mathcal{T}_X}  \overline{ \mathbf E}_x \! \left[1_{\{\tau<\tau_{X \geq Y} \! \}}\, \mathrm e^{-r \tau} R(X_{\tau}) +1_{\{\tau \geq \tau_{X \geq Y} \! \}}\, \mathrm e^{-r \tau_{X \geq Y}}\! \left[V_U(X_{\tau_{X\geq Y}})-P(X_{\tau_{X \geq Y}})\right] \right] \hskip -1mm.
\end{eqnarray*}
The proof consists of two parts.

\subparagraph{Proof of $\geq$}

Consider the filtrations $(\mathcal F^0_t) _{t\geq 0} $ and $(\mathcal F_t) _{t\geq 0} $ over $\Omega$ defined by $\mathcal{F}^0_t \equiv \sigma(X_s; s\leq t)$ and $\mathcal{F}_t \equiv \cap_{s> t}\mathcal{F}^0_s$ for all $t \geq 0$ and let $\mathcal{T}^0_X$ denote the set of stopping times with respect to the filtration $(\mathcal F^0_t) _{t\geq 0} $. Recall that the filtration $(\mathcal{G}_t)_{t \geq 0}$ over $\overline \Omega$ is defined by $\mathcal{G}_t \equiv \sigma( X_s,  1_{\{\tau_{X \geq Y} \leq s\}}; s\leq t)$ for all $t \geq 0$, so that we have, for any such $t$,
\begin{eqnarray}
\mathcal{F}^0_t\otimes \{\emptyset, \mathcal I\} \subset \mathcal{G}_t  \subset \mathcal{F}^0_t \otimes \mathcal{B}(I) \label{betweennessmathcalG}
\end{eqnarray}
and
\begin{eqnarray}
\mathcal{G}_t = \sigma \!\left(X_s, \{Y\leq X_0\}, \left\{Y >  \sup_{s \in [0,t]} X_s\right\} \!, Y 1_{\{Y \in (X_0, \,\sup_{s \in [0,t]}X_s]\}};  s \leq t\right) \hskip -1mm , \label{mathcalG}
\end{eqnarray}
and $\tau_{X \geq Y}$ is a stopping time with respect to the filtration $(\mathcal{G}_t)_{t \geq 0}$. Notice also that $\mathcal{T}^0_X \subset \mathcal{T}_{X, X \geq Y}$ if we identify the elements of $\mathcal{T}^0_X$ to functions defined on $\overline \Omega$. For any $\tau \in \mathcal{T}_{X,X \geq Y}$, let us define
\begin{eqnarray*}
\tilde \tau (\tau) \equiv \tau 1_{\{\tau < \tau_{X \geq Y}\! \}}+ (\tau_{X\geq x_U}\vee \tau_{X \geq Y})  1_{\{\tau \geq \tau_{X \geq Y}\!\}}.
\end{eqnarray*}
Using the properties of the $\sigma$-fields $\mathcal{G}_\tau$ and $\mathcal{G}_{\tau_{X \geq Y}}$, we have $\tilde \tau(\tau) \in \mathcal{T}_{X, X \geq Y}$ as
\begin{eqnarray*}
\{\tilde \tau (\tau) \leq t \}= (\{ \tau \leq t \}\cap  \{\tau < \tau_{X \geq Y}\}  ) \cup (\{ \tau_{X\geq x_U} \leq t \}\cap \{\tau_{X \geq Y} \leq t\} \cap  \{\tau \geq \tau_{X \geq Y}\})
\end{eqnarray*}
for all $t\geq 0$. Notice that
\begin{eqnarray*}
\tau_{X\geq x_U}\vee \tau_{X \geq Y} = \tau_{X\geq Y} + \tau_{X \geq x_U} \circ \theta_{\tau_{X \geq Y}} 1_{\{\tau_{X \geq Y}<\infty\}},
\end{eqnarray*}
where $\theta_.$ denotes the shift operator on $\Omega$.

By construction, $X$ is a strong Markov process with respect to the filtration generated by $X$, and also with respect to the filtration generated by $X$ and $Y$ as $Y$ is independent of $X$. Hence, because the filtration $(\mathcal{G}_t)_{t \geq 0}$ lies in between these two filtrations by \eqref{betweennessmathcalG}, $X$ is also a strong Markov process with respect to the filtration $(\mathcal{G}_t)_{t \geq 0}$. Denoting by $\tilde X$ another copy of the canonical process defined on $(\tilde \Omega,\tilde{ \mathcal{F}} )=(\Omega,\mathcal {F})$, the strong Markov property yields
\begin{align*}
\hat J (x ,& \tilde\tau (\tau) )
\\
&=\overline{ \mathbf E}_x \! \left[1_{\{\tau<\tau_{X \geq Y} \! \}} \, \mathrm e^{-r \tau} R(X_{\tau}) \right.
\\
&   \; \; \; \; \; \;\; \;\,\; \left. +\,1_{\{\tau \geq\tau_{X \geq Y}\! \}}\, \mathrm e^{-r \tau_{X \geq Y}} \! \left[ \mathrm e^{-r\tau_{X \geq x_U} \circ \,\theta_{\tau_{X \geq Y}}} U(X_{\tau_ {X\geq Y}+\,\tau_{X \geq x_U} \circ \,\theta_{\tau_{X \geq Y}}} ) - P(X_{\tau_{X \geq Y}}) \right] \right]
\\
&=  \overline{ \mathbf E}  _x \! \left[1_{\{\tau<\tau_{X \geq Y} \! \}}\, \mathrm e^{-r \tau} R(X_{\tau})+1_{\{ \tau \geq \tau_{X \geq Y}\! \}} \, \mathrm  e^{-r \tau_{X \geq Y}} \, \mathbf {E}_{X_{\tau_{X \geq Y}}} \! \! \left[\mathrm e^{-r\tau_{\tilde X \geq x_U}} U(\tilde X_{\tau_{\tilde X \geq x_U}} ) - P(\tilde X_0) \right] \right]
\\
&= \overline{ \mathbf E}  _x \! \left[1_{\{\tau<\tau_{X \geq Y} \! \}}\, \mathrm e^{-r \tau} R(X_{\tau})+1_{\{\tau \geq\tau_{X \geq Y} \! \}} \, \mathrm e^{-r \tau_{X \geq Y}}\!   \left[V_U(X_{\tau_{X\geq Y}})-P(X_{\tau_{X\geq Y}})\right] \right] \hskip -1mm .
\end{align*}
We deduce from this that
\begin{align*}
 \sup_{\tau \in \mathcal{T}_{X, X\geq Y}} & \hat J(x, \tau)
\\
& \geq \sup_{\tau \in \mathcal{T}_{X, X\geq Y}} \hat J(x, \tilde \tau(\tau))
\\
&=   \sup_{\tau \in \mathcal{T}_{X, X\geq Y}}\overline{ \mathbf E}  _x \! \left[1_{\{\tau<\tau_{X \geq Y} \! \}}\, \mathrm e^{-r \tau} R(X_{\tau }) +1_{\{\tau \geq\tau_{X \geq Y} \! \}} \, \mathrm e^{-r \tau_{X \geq Y}}\!  \left[V_U(X_{\tau_{X\geq Y}})-P(X_{\tau_{X\geq Y}})\right] \right]
\\
& \geq  \sup_{\tau \in \mathcal{T}^0_X}\overline{ \mathbf E}  _x \! \left[1_{\{\tau<\tau_{X \geq Y} \! \}}\, \mathrm e^{-r \tau} R(X_{\tau})+1_{\{\tau \geq\tau_{X \geq Y} \! \}} \, \mathrm e^{-r \tau_{X \geq Y}}\!\left[V_U(X_{\tau_{X\geq Y}})-P(X_{\tau_{X\geq Y}})\right] \right]
\\
&=  \sup_{\tau \in \mathcal{T}_X}\overline{ \mathbf E}  _x \! \left[1_{\{\tau<\tau_{X \geq Y} \! \}}\, \mathrm e^{-r \tau} R(X_{\tau})+1_{\{\tau \geq\tau_{X \geq Y} \! \}} \, \mathrm e^{-r \tau_{X \geq Y}}\!\left[V_U(X_{\tau_{X\geq Y}})-P(X_{\tau_{X\geq Y}})\right] \right] \hskip -1mm,
\end{align*}
where the first inequality follows from the fact that $\tilde \tau(\tau) \in \mathcal{T}_{X, X \geq Y}$ for all $ \tau \in \mathcal{T}_{X, X \geq Y}$, the second inequality follows from $\mathcal{T}^0_X \subset \mathcal{T}_{X, X\geq Y}$, and the second equality follows from the fact that any $\tau \in \mathcal{T}_X$ is the limit of the nonincreasing sequence of stopping times $(\tau_n) _{n \in \mathbb N}$ in $\mathcal{T}^0_X$ defined by
\begin{eqnarray*}
\tau_n \equiv \sum_{k \geq 0} 1_{\{\tau \in [k2^{-n}, (k+1)2^{-n})\}} (k+1)2^{-n}, \quad n \in \mathbb N,
\end{eqnarray*}
which, given that $X$ is continuous and $R$ is continuous and satisfies A1, allows us to apply Lebesgue's dominated convergence theorem to replace the supremum over $\mathcal{T}^0_X$ with the supremum over $\mathcal T_X$. This concludes the first part of the proof.

\subparagraph{Proof of $\leq$}

The proof of the reverse inequality is more technical, although it relies on a very intuitive decomposition of stopping times in $\mathcal{T}_{X, X \geq Y}$. Specifically, we show that for any $ \tau \in \mathcal{T}_{X, X \geq Y}$, there exists $\tau^0 \in \mathcal{T}^0_X$ such that $\{\tau< \tau_{X \geq Y}\}=\{\tau^0< \tau_{X \geq Y}\}$ and
\begin{eqnarray} \label{stopping-decomposition}
\tau(\omega,y)= \tau^0(\omega) 1_{\{\tau^0(\omega)< \tau_{X \geq Y}(\omega, y)\}} + \tau(\omega,y) 1_{\{\tau^0 (\omega) \geq \tau_{X \geq Y}(\omega,y)\}}, \quad  (\omega,y)\in \overline \Omega .
\end{eqnarray}
We provide a detailed proof for lack of a reference covering exactly the case at hand. The proof is constructive. Let us choose a continuous strictly increasing map $\psi: \mathcal I \rightarrow \mathcal I$ such that $\psi(x)>x$ for all $x\in \mathcal I$, and define
\begin{eqnarray*}
\phi(\omega,t)\equiv\tau \! \left(\omega,\psi \!\left(\sup_{s \in [0,t]}\omega_s\right)\right)\hskip -1mm, \quad (\omega,t) \in \Omega \times \mathbb R_+.
\end{eqnarray*}
We verify that \eqref{stopping-decomposition} holds for
\begin{eqnarray*}
\tau^0(\omega) \equiv \inf \hskip 0.5mm \{ t \geq 0 : \phi(\omega,t)\leq t \}, \quad \omega \in \Omega.
\end{eqnarray*}
To see this, observe that for each $t\in \mathbb{R}_+$, the mapping $(\omega,y) \mapsto 1_{\{\tau(\omega,y) \leq t\}}$ is $\mathcal{G}_t$-measurable. By (\ref{mathcalG}), this implies that for each $(\omega,t) \in \Omega \times \mathbb{R}_+$, the mapping $y \mapsto 1_{\{\tau(\omega,y)\leq t \}}$ is constant over $(\sup_{s \in [0,t]}\omega_s, \beta)$. Because $\psi(x) > x$ for all $x \in \mathcal I$, it follows that, for each $t\geq 0$, the following equivalences hold:
\begin{align}
\phi(\omega,t)\leq t  &\Leftrightarrow  \mbox{there exists $y> \sup_{s \in [0,t]}\omega_s$ such that $\tau(\omega,y) \leq t$}  \notag
\\
 &\Leftrightarrow  \mbox{for each $y> \sup_{s \in [0,t]}\omega_s$, $\tau(\omega,y) \leq t$}. \label{propertyG}
\end{align}
A consequence of this is that, for each $\omega$ such that $\tau^0(\omega)<\infty$,
\begin{eqnarray} \label{tau0-interval}
\{ t \geq 0 : \phi(\omega,t)\leq t \}=[\tau^0(\omega),\infty).
\end{eqnarray}
Indeed, assume that $\tau^0(\omega)<\infty$. For each $t' \geq t$, $\psi(\sup_{s \in[0, t']} \omega _s) >  \sup_{s \in[0, t']} \omega _s \geq \sup_{s \in[0, t]} \omega _s$. Hence, by \eqref{propertyG}, $\phi(\omega,t)\leq t$ implies that, for each $t' \geq t$, $\phi(\omega,t')\leq t\leq t'$, and thus $(\tau^0(\omega),\infty) \subset \{ t \geq 0 : \phi(\omega,t)\leq t \}$. Finally, if $y>\sup_{s \in [0,\tau^0(\omega)]}\omega_s$, then $y>\sup_{s \in [0,t']}\omega_s$ for $t'>\tau^0(\omega)$ sufficiently close to $\tau^0( \omega)$, and thus $\tau(\omega,y) \leq t'$ by \eqref{propertyG} as $\phi(\omega,t')\leq t'$. This implies $\tau(\omega,y) \leq \tau^0(\omega)$ by taking the limit and thus $\phi(\omega,\tau^0(\omega))\leq \tau^0(\omega)$ by \eqref{propertyG} again, which concludes the proof of \eqref{tau0-interval}. Now, if $\tau_{X \geq y}(\omega) > \tau^0(\omega)$, then $y>\sup_{s \in [0,\tau^0(\omega)]}\omega_s$ and, as above, we obtain that $\tau(\omega,y)\leq \tau^0(\omega)$. Similarly, if $\tau_{X \geq y}(\omega) > \tau(\omega,y)$, then $y>\sup_{s \in [0,\tau(\omega,y)]}\omega_s$,  implying that  $\phi(\omega,\tau(\omega,y))\leq \tau(\omega,y)$ by \eqref{propertyG}, and thus that $\tau^0(\omega)\leq \tau(\omega,y)$ by (\ref{tau0-interval}).
These two inequalities together imply
\begin{eqnarray}\label{equality_before_threshold}
\{\tau_{X \geq y} > \tau\} = \{\tau_{X \geq y} > \tau^0\}  \text{ and } \tau(\omega,y)=\tau^0(\omega) \text{ if } \tau_{X \geq y}(\omega) > \tau(\omega,y) .
\end{eqnarray}
For each $t \in \mathbb R_+$, we have $\{ \tau^0 \leq t\}= \{ \phi(\omega,t) \leq t\}$ by \eqref{tau0-interval}, and this set belongs to $\mathcal{F}^0_t$ because  $\{ (\omega,y) : \tau(\omega,y) \leq t\} \in \mathcal{F}^0_t \otimes \mathcal{B}(\mathcal I)$ and the mapping $\omega \mapsto (\omega,\psi(\sup_{[0,t]}\omega_s))$ is $\mathcal{F}^0_t/\mathcal{F}^0_t \otimes \mathcal{B}(I)$ measurable. It follows that $\tau^0 \in \mathcal{T}^0_X$, which, together with \eqref{equality_before_threshold}, concludes the proof of \eqref{stopping-decomposition}.

Using this decomposition, we have
\begin{align*}
\hat J( x , \tau) &= \overline{ \mathbf E}  _x  \! \left[1_{\{\tau<\tau_{X \geq Y} \! \}}\, \mathrm e^{-r \tau} R(X_{\tau}) +1_{\{\tau \geq \tau_{X \geq Y}\! \}} \! \left[ \mathrm e^{-r \tau} U(X_{\tau} ) - \mathrm e^{-r\tau_{X \geq Y}}P(X_{\tau_{X \geq Y}}) \right] \right]
\\
&=\overline{ \mathbf E}  _x \! \left[1_{\{\tau^0<\tau_{X \geq Y} \! \}} \, \mathrm e^{-r \tau^0} R(X_{\tau^0}) +\overline{ \mathbf E}_x \!\left[ 1_{\{\tau^0 \geq \tau_{X \geq Y\!}\}} \! \left[ \mathrm e^{-r \tau} U(X_{\tau}) - \mathrm e^{-r\tau_{X \geq Y}}P(X_{\tau_{X \geq Y}}) \right) \! \mid \!Y\right ] \right] \hskip -1mm.
\end{align*}
Notice that, by definition, $\tau_{X \geq Y}(\omega,y)=\tau_{X \geq y}(\omega)$ and thus $\tau_{X \geq y} \in \mathcal{T}^0_X$ for all $y \in \mathcal I$. Moreover, for any such $y$, $\tau_y(\cdot)\equiv \tau(\cdot,y) \in \mathcal{T}^0_X$ because $\tau$ is also a stopping time with respect to the larger filtration $(\mathcal{F}^0_t\otimes \mathcal{B}(\mathcal I))_{t \geq 0}$. It follows from Dellacherie and Meyer (1975, Th\'eor\`eme 103) that, for each $y \in \mathcal I$, there exists a $\mathcal{F}^0_{\tau_{X \geq y}}\otimes \mathcal{F}^0_\infty$-measurable function $T_y:\Omega \times \Omega \rightarrow [0,\infty]$ such that $T_y(\omega, \cdot) \in \mathcal{T}^0_X$ for all $\omega \in \Omega$ and
\begin{eqnarray*}
\tau_y(\omega)1_{\{ \tau_y (\omega) \geq \tau_{X \geq y}(\omega) \}} = [\tau_{X \geq y} (\omega)+T_y(\omega, \theta_{\tau_{X \geq y}}(\omega))] 1_{\{ \tau_y (\omega) \geq \tau_{X \geq y} (\omega) \}}
\end{eqnarray*}
for all $\omega \in \Omega$ such that $\tau_{X \geq y}(\omega)<\infty$. We obtain that, for each $y \in \mathcal I$,
\begin{align*}
\overline{ \mathbf E}_x \! \left[ 1_{\{\tau^0 \geq \tau_{X \geq Y} \}} \, \mathrm e^{-r \tau} U(X_{\tau} ) \! \mid \! Y=y \right] \! &= \mathbf E_x \! \left[ 1_{\{\tau^0 \geq \tau_{X \geq y} \}}\, \mathrm e^{-r \tau_y} U(X_{\tau_y} )\right]
\\
&=\mathbf E_x \! \left[ 1_{\{\tau^0 \geq \tau_{X \geq y} \}} \, \mathrm e^{-r ( \tau_{X \geq y} + T_y)} U(X_{\tau_{X \geq y} +T_y} )\right]
\\
&=\mathbf E_x \! \left[ 1_{\{\tau^0 \geq \tau_{X \geq y} \}} \, \mathrm e^{-r \tau_{X \geq y}} \,\mathbf E_x \! \left[\mathrm e^{-r T_y} U(X_{\tau_{X \geq y} +T_y} ) \! \mid \! \mathcal{F}^X_{\tau_{X \geq y}}\right]\right]\hskip -1mm.
\end{align*}
Applying the Markov property yields
\begin{align*}
\mathbf E_x \! \left[ \mathrm e^{-r T_y} U(X_{\tau_{X \geq y} +T_y} ) \! \mid \! \mathcal{F}^X_{\tau_{X \geq y}}\right]\! (\omega) &=\int_{\Omega} \mathrm e^{-r T_y(\omega,\tilde \omega)} U(X_{T_y(\omega,\tilde\omega)}) \, \mathbf P_{X_{\tau_{X \geq y}}(\omega)}(\mathrm d\tilde\omega)
\\
& \leq V_U(X_{\tau_{X \geq y}}(\omega))
\end{align*}
for all $(\omega, y) \in \Omega \times \mathcal I$ such that $\tau_{X \geq y}(\omega)<\infty$. Summing up, we obtain
\begin{align*}
\hat J(x, \tau) &\leq \overline{ \mathbf E}  _x \! \left[1_{\{\tau^0<\tau_{X \geq Y} \! \}}\, \mathrm e^{-r \tau^0} R(X_{\tau^0}) +1_{\{ \tau^0 \geq \tau_{X \geq Y} \! \}}\, \mathrm e^{-r\tau_{X \geq Y}} \!\left[V_U(X_{\tau_{X \geq Y}})-P(X_{\tau_{X \geq Y}})\right] \right] \\
&\leq \sup_{\tilde \tau \in \mathcal{T}_X} \overline{ \mathbf E}  _x \! \left[1_{\{\tilde\tau<\tau_{X \geq Y} \! \}} \, \mathrm e^{-r \tilde\tau} R(X_{\tilde\tau}) +1_{\{ \tilde \tau \geq \tau_{X \geq Y} \! \}}\, \mathrm  e^{-r\tau_{X \geq Y}} \!\left[V_U(X_{\tau_{X \geq Y}})-P(X_{\tau_{X \geq Y}}) \right]  \right]\hskip -1mm,
\end{align*}
where the second inequality follows from $\tau^0 \in \mathcal T_X^0 \subset \mathcal T_X$. This concludes the second part of the proof. The result follows. \hfill $\blacksquare$

\bigskip

It should be noted that the decomposition (\ref{stopping-decomposition}) of stopping times in $\mathcal T_{X, X \geq Y}$ allows us to decompose any solution to (\ref{11-2})---where the supremum is taken over the stopping times in $\tau_{X, X \geq Y}$---into a solution to \eqref{mainpb-intro}---where the supremum is taken over the stopping times in $\tau_X$---and a solution to (\ref{saU}) in the continuation game.

\renewcommand{\thesection}{Appendix C: An Example}

\section{}

\renewcommand{\thesection}{C}

\renewcommand{\theequation}{C.\arabic{equation}}
\setcounter{equation}{0}

In this appendix, we verify that the specification of the model provided at the end of Section \ref{TecAss} satisfies A1--A8. We let $X$ follow a geometric Brownian motion with drift $\mu<r$ and volatility $\sigma >0$,
\begin{eqnarray*}
\mathrm dX_t= \mu X_t \, \mathrm dt + \sigma X_t \, \mathrm dW_t, \quad t\geq 0,
\end{eqnarray*}
so that the state space for $X$ is $\mathcal I= (0, \infty)$ and the infinitesimal generator of $X$ writes as
\begin{eqnarray*}
\mathcal L u(x) \equiv \mu x u'(x) + \frac{1}{2}\, \sigma^2 x^2 u''(x), \quad x \in (0, \infty).
\end{eqnarray*}
The two fundamental solutions to $\mathcal L u -ru = 0$ are, up to a linear transformation,
\begin{eqnarray*}
h_1 (x) \equiv x^{\beta_1} \quad \mbox{and} \quad h_2(x) \equiv x ^{\beta_2}, \quad x \in (0, \infty),
\end{eqnarray*}
where, letting  $\nu \equiv \frac{\mu}{\sigma^2}  - \frac{1}{2}$,
\begin{eqnarray*}
\beta_1 \equiv -  \nu + \sqrt{\nu^2 + \frac{2 r }{\sigma^2} } >1 \quad \mbox{and} \quad \beta_2 \equiv - \nu - \sqrt{\nu^2 + \frac{2 r }{\sigma^2} } <0.
\end{eqnarray*}
Observe that $h_1$ and $h_2$ are strictly convex over $(0, \infty)$. The derivative of the scale function is, up to a linear transformation,
\begin{eqnarray*}
S'(x) = \exp \! \left ( - \int_1^x \frac{2 \mu}{\sigma^2 z}  \, \mathrm dz \right ) \! = x^{-2 \nu - 1}, \quad x \in (0,\infty),
\end{eqnarray*}
which leads to
\begin{eqnarray*}
\gamma = \frac{h_1'(x) h_2(x)  - h_1(x)  h_2'(x)}{S'(x)} = \beta_1 - \beta_2.
\end{eqnarray*}
The payoff functions are $R(x) \equiv x-I$ and $U(x) \equiv \kappa x -I$ for $\kappa >1$, which satisfy A3 for $x_0 ={r \over r- \mu}\, I$ and $x_0 = {r \over r - \mu}\, I \kappa ^{-1}$, respectively. That $R$ and $U$ satisfy A2 follows from the explicit expression
\begin{eqnarray*}
X_t = x \exp \!\left(\!\left( \mu-{1\over 2} \, \sigma^2 \right)\!t + \sigma W_t\right)
\end{eqnarray*}
for the geometric Brownian motion starting at $X_0 =x$ along with the assumption that $r > \mu$. That $R$ and $U$ satisfy A1 follows from observing that
\begin{align*}
\mathbf E _x \! \left[ \sup_{t \geq 0} \,\mathrm e^{-rt} X_t \right]  \! &=x \, \mathbf E_x \! \left[ \left[\exp \! \left( \sup_{t \geq 0} \, W_t - {1 \over \sigma} \!\left( r- \mu +{1\over 2}\,  \sigma^2 \right) \!t  \right) \right]^\sigma \right]
\\
& =   {2 \over \sigma} \! \left( r- \mu +{1\over 2}\,  \sigma^2 \right) x \int _0 ^\infty \mathrm \exp \! \left( \sigma y -   {2 \over \sigma} \!\left( r- \mu +{1\over 2}\,  \sigma^2 \right) \! y \right) \mathrm dy
\\
& =   {r- \mu +{1\over 2}\,  \sigma^2 \over r- \mu}\, x,
\end{align*}
where the second inequality follows from the fact that, for each $\lambda >0$, the random variable $\sup _{t \geq 0} W_t-\lambda t$ has an exponential density with parameter $2\lambda$ (Revuz and Yor (1999, Chapter II, \S3, Exercise 3.12)). The value functions $V_R$ and $V_U$ write as
\begin{eqnarray*}
V_R(x) = \left\{ \begin{array}{lll} \frac{h_1(x)}{h_1(x_R)}  \,(x_R -I)& \text{if} & x < x_R, \\  x - I & \text{if} & x \geq x_R, \end{array} \right. \quad \mbox{and} \quad V_U(x) = \left\{ \begin{array}{lll} \frac{h_1(x)}{h_1(x_U)} \,(\kappa x_U -I )& \text{if} & x < x_U, \\  \kappa x - I & \text{if} & x \geq x_U, \end{array} \right.
\end{eqnarray*}
where
\begin{eqnarray*}
x_R \equiv \frac{\beta_1}{\beta_1 - 1} \, I \quad \mbox{and} \quad x_U \equiv \frac{\beta_1}{\beta_1 - 1}\, I \kappa^{-1}.
\end{eqnarray*}
Observe that $x_R > x_U$ and that $V_R$ and $V_U$ are $\mathcal C^1$ over $(0, \infty)$, that is, $V_R$ and $V_U$ satisfy the smooth-fit property. Moreover, $V_R$ and $V_U$ are $\mathcal C^2$ and satisfy $\mathcal LV_R - rV_R \leq 0$ and $\mathcal LV_U - rV_U \leq 0$ over $(0, \infty) \setminus \{x _R \}$ and $(0, \infty) \setminus \{x_U \}$, respectively, so that the function $G = \frac{1}{2} \, (V_U+ V_R)$  is $\mathcal C^2$ and satisfies A7 over $(0, \infty) \setminus \{x_R, x_U\}$. That $G$ satisfies A5 follows again from the explicit expression for the geometric Brownian motion along with the fact that $G$ is bounded above by a linear \pagebreak function, and that $G$ satisfies A6 follows from the fact that $V_U > V_R$. Finally, because
\begin{eqnarray*}
0 < G   < {1 \over 2} \,[|R | +R(x_R)  + |U | + U(x_U)]
\end{eqnarray*}
over $(0, \infty)$, $G$ satisfies A4 because $R$ and $U$ satisfy A1. Now, consider the function $P = \frac{1}{2}\, (V_U   - V_R ) $. According to the analysis in Section \ref{AMotExa}, we need to show that $P$ is strictly increasing and onto. We have
\begin{eqnarray}  \label{K}
V_U(x) -V_R(x) = \left\{ \begin{array}{lll} \frac{h_1(x)}{h_1(x_U)}\,(\kappa x_U -I ) - \frac{h_1(x)}{h_1(x_R)}\, (x_R -I)& \text{if} & x < x_U, \\
\kappa x -I -\frac{h_1(x)}{h_1(x_R)}\, (x_R -I) & \text{if} &  x_U\leq x < x_R,\\
(\kappa -1) x & \text{if} & x \geq x_R. \end{array} \right.
\end{eqnarray}
That $P$ is strictly increasing over $[x_R, \infty)$ is obvious. Next, because $V_U > V_R$, we have
\begin{eqnarray*}
\frac{1}{h_1(x_U)}\,(\kappa x_U -I ) > \frac{1}{h_1(x_R)}\, (x_R -I),
\end{eqnarray*}
which implies, as $h_1$ is strictly increasing, that $P$ is strictly increasing over $(0, x_U)$. Finally, for each  $x \in [x_U, x_R)$, we have
\begin{eqnarray*}
V'_U(x)-V'_R(x) = \kappa - \frac{h'_1(x)}{h_1(x_R)}\, (x_R -I) >1 - \frac{h'_1(x_R)}{h_1(x_R)}\, (x_R -I) =0
\end{eqnarray*}
where the inequality follows from the fact that $h_1$ is strictly convex and that $x_R > I$, and the second equality follows from the smooth-fit property for $V_R$. This implies that $P$ is strictly increasing over $[x_U, x_R)$. Because $P$ is continuous, and because $\lim_{x \to 0^+}P(x) = 0$ and $\lim_{x \to \infty}P(x) = \infty$, we obtain that $P$ is strictly increasing and onto, as desired. To conclude, observe that, if $Z$ is drawn from a distribution with locally Lispchitz density $f_Z>0$ over $(0, \infty)$ with respect to Lebesgue measure, then the density of the law $\mathbf Q$ of $Y = P^{-1}(Z)$ is given by $f =  P'  f_Z \circ P >0$. It then follows from (\ref{K}) and from the properties of the value functions $V_R$ and $V_U$ that $P' \in \mathcal C^1((0, \infty) \setminus \{ x_R, x_U \})$, with bounded left- and right-derivatives at $x_R$ and $x_U$, so that A8 is satisfied.

\end{document}